\documentclass[11pt,a4paper]{amsart}

\usepackage[utf8]{inputenc}
\usepackage[T1]{fontenc}
\usepackage[activeacute,english]{babel}
\usepackage{lmodern}
\usepackage{microtype}
\usepackage{marvosym}
\usepackage{stmaryrd}
\usepackage{amsaddr}

\usepackage{amsmath,amssymb,amsthm,amsfonts,mathtools}

\usepackage{enumitem}
\usepackage[colorlinks,bookmarks,linkcolor=black,citecolor=black]{hyperref}

\usepackage{xspace,caption,subcaption,comment}
\usepackage{multirow}
\usepackage{xcolor}
\usepackage[normalem]{ulem}
\usepackage{fullpage}

\newtheorem{theorem}{Theorem}[section]
\newtheorem{lemma}[theorem]{Lemma}
 
\newtheorem{conj}[theorem]{Conjecture}
\newtheorem{question}[theorem]{Question}

\newtheorem{claim}[theorem]{Claim}

\newtheorem{definition}[theorem]{Definition}

\DeclareMathOperator{\im}{Im}
\DeclareMathOperator{\ORB}{ORB}
\DeclareMathOperator{\orb}{Orb}
\DeclareMathOperator{\stab}{Stab}

\makeatletter

\newif\ifdraft
\draftfalse
\ifdraft 
    \usepackage[notcite,notref]{showkeys}

  \else
    \renewcommand{\sout}[1]{} %making command \sout hide the text in the final version: in the draft version striking out is to emphasize a recent change
\fi

\newcommand{\R}{\mathbb{R}}
\newcommand{\N}{\mathbb{N}}
\newcommand{\eps}{\varepsilon}

\newcommand{\nth}{\varnothing}
\newcommand{\sm}{\setminus}
\newcommand{\ca}{\mathcal{A}}
\newcommand{\cb}{\mathcal{B}}
\newcommand{\cc}{\mathcal{C}}
\newcommand{\cd}{\mathcal{D}}

\newcommand{\cg}{\mathcal{G}}
\newcommand{\ch}{\mathcal{H}}
\newcommand{\ci}{\mathcal{I}}
\newcommand{\cj}{\mathcal{J}}
\newcommand{\cp}{\mathcal{P}}
\newcommand{\cq}{\mathcal{Q}}
\newcommand{\cs}{\mathcal{S}}

\newcommand{\D}{\,\mathrm{d}}
\newcommand{\rul}{\mathcal{R}}
\newcommand{\lgr}[1]{\ch_{#1}}
\newcommand{\Kernel}{\mathcal{W}}
\newcommand{\Gra}{\Kernel_{0}}
\newcommand{\indic}{\mathbb{I}}
\newcommand{\ind}[1]{\indic_{ \left\{ #1 \right\}}}

\newcommand{\Troot}[1]{T_{#1}}
\newcommand{\vel}[1][\:]{\mathfrak{V}_{#1}}
\newcommand{\traj}[2]{\Phi^{#1}#2}
\newcommand{\trajrul}[3]{\Phi_{#1}^{#2}#3}
\newcommand{\mdom}[1]{\cd_{#1}}
\newcommand{\life}[1]{L_{#1}}
\newcommand{\tf}{\mathrm{tf}}
\newcommand{\rt}{\mathrm{rt}}
\newcommand{\m}{\mathbf{m}}
\newcommand{\n}{\mathbf{n}}
\newcommand{\p}{\mathbf{p}}
\newcommand{\q}{\mathbf{q}}
\newcommand{\vv}{\mathbf{v}}
\newcommand{\x}{\mathbf{x}}
\newcommand{\y}{\mathbf{y}}
\newcommand{\z}{\mathbf{z}}
\newcommand{\A}{\mathbf{A}}
\newcommand{\F}{\mathbf{F}}
\newcommand{\G}{\mathbf{G}}
\newcommand{\HH}{\mathbf{H}}
\newcommand{\PP}{\mathbf{P}}

\newcommand{\Om}{\mathbf{\Omega}}

\newcommand{\cutnd}{d_\square}

\newcommand{\Linf}[1]{ \left\lVert #1\right\rVert_{\infty}}

\def\endofClaim{\hfill\scalebox{.6}{$\blacksquare$}}
\newcommand{\oldqed}{}
\newenvironment{claimproof}[1][Proof]{
  \renewcommand{\oldqed}{\qedsymbol}
  \renewcommand{\qedsymbol}{\endofClaim}
  \begin{proof}[#1]
}{
  \end{proof}
  \renewcommand{\qedsymbol}{\oldqed}
}

\title{Characterization of flip process rules with the same trajectories}
\author{Eng Keat Hng}
\email{hng@ibs.re.kr}
\address{Extremal Combinatorics and Probability Group, Institute for Basic Science, 55 Expo-ro, Yuseong-gu, Daejeon 34126, South Korea}
\thanks{Supported by IBS-R029-C4. Most of the work was done while affiliated with the Institute of Computer Science of the Czech Academy of Sciences (institutional support RVO:67985807) and supported by Czech Science Foundation Project 21-21762X.}

\begin{document}

\begin{abstract}
Garbe, Hladk\'{y}, \v{S}ileikis and Skerman [Ann.\ Inst.\ Henri Poincaré Probab.\ Stat., 60 (2024), pp.\ 2878--2922] recently introduced a general class of random graph processes called flip processes and proved that the typical evolution of these discrete-time random graph processes corresponds to certain continuous-time deterministic graphon trajectories. We obtain a complete characterization of the equivalence classes of flip process rules with the same graphon trajectories. As an application, we characterize the flip process rules which are unique in their equivalence classes. These include several natural families of rules such as the complementing rules, the component completion rules, the extremist rules, and the clique removal rules.
\end{abstract}

\maketitle

\section{Introduction}

Graphs are mathematical structures which underpin the study of numerous real-life settings, many of which evolve in time according to a stochastic rule. This leads us to \emph{random graph processes}, the study of which was initiated by Erd\H{o}s and R\'{e}nyi in their seminal paper~\cite{ErdosRenyi60} from 1960. In the decades that followed, a variety of random graph processes were introduced and studied extensively. One example is the triangle removal process, which was introduced by Bollob\'{a}s and Erd\H{o}s~\cite{Bollobas98} in 1990. The process begins with the initial graph $G_0 = K_n$, and in each step $\ell$ we obtain the graph $G_{\ell}$ by deleting the edges of a uniformly random triangle in $G_{\ell-1}$. The process terminates when $G_{\ell}$ is triangle-free. There is a long line of research on random graph processes resulting in triangle-free graphs, which has led to a series of breakthroughs on lower bounds for the Ramsey numbers $R(3,t)$ in~\cite{Kim,Bohman,BohmanKeevash,FizPontiverosGriffithsMorris}.

Consider the following modification of the triangle removal process. Instead of deleting the edges of a uniformly random triangle, in each step we sample a uniformly random triple of distinct vertices, delete the edges of any triangle induced and do nothing otherwise. Note that this has the form of a `local replacement' rule: triples of distinct vertices represent `localities' where we replace triangles by empty graphs and do nothing otherwise. Furthermore, observe that we may recover the original triangle removal process from the modified process by ignoring the steps that do nothing. Indeed, conditioned on sampling a triangle, every step of the modified process deletes a uniformly random triangle.

Generalising the idea of local replacements, Garbe, Hladk\'{y}, \v{S}ileikis and Skerman~\cite{GarbeHladkySileikisSkerman} recently introduced a general class of random graph processes called flip processes. Before we describe these processes, we first introduce some notation and terminology. Write $\N$ for the set of positive integers and $\N_0$ for the set $\N\cup\{0\}$. For $a\in\N_0$ write $[a]$ for the set $\{1,\dots,a\}$ and $[a]_0$ for $[a]\cup\{0\}$. All our graphs are simple and undirected. Let $k\in\N$ and let $\lgr{k}$ be the set of all labelled graphs on the vertex set $[k]$. A \emph{rule} is a matrix $\rul = (\rul_{F,H})_{F,H\in\lgr{k}} \in [0,1]^{\lgr{k}\times\lgr{k}}$ such that for each $F\in\lgr{k}$ we have $\sum_{H\in\lgr{k}} \rul_{F,H} = 1$. We call $k$ the \emph{order} of the rule $\rul$. A \emph{flip process} with rule $\rul$ is a random graph process $(G_{\ell})_{\ell\in\N_0}$  where we start with an \emph{initial graph} $G_0$ on the vertex set $[n]$ with $n \ge k$ and in step $\ell\in\N$ we obtain the graph $G_{\ell}$ from $G_{\ell-1}$ as follows. First, we sample an ordered tuple $\vv = (v_i)_{i\in[k]}$ of distinct vertices in $[n]$ at random. Next, we define the \emph{drawn graph} $F\in\lgr{k}$ by $ij \in E(F)$ if and only if $v_iv_j\in E(G_{\ell-1})$ and generate the \emph{replacement graph} $H\in\lgr{k}$ according to the probability distribution $(\rul_{F,H})_{H\in\lgr{k}}$. Finally, we replace $G_{\ell-1}[\vv]$ with $H$ to obtain $G_{\ell}$ from $G_{\ell-1}$.

Flip processes were studied in~\cite{GarbeHladkySileikisSkerman} through the lens of dense graph limits~\cite{LovaszBook}. Let $(\Omega,\pi)$ be an atomless probability space with an implicit separable sigma-algebra. Write $\Kernel \subseteq L^\infty(\Omega^2)$ for the set of \emph{kernels}, that is, bounded symmetric measurable real-valued functions. Write $\Gra$ for the set of \emph{graphons}, that is, the elements of $\Kernel$ whose range is a subset of $[0,1]$. An important metric on $\Kernel$ is the \emph{cut norm distance} $\cutnd$ given by $\cutnd(U,W) = \sup_{S,T\subseteq\Omega} \left|\int_{S\times T}(U-W)\D\pi^2\right|$. Let $G$ be a graph and $(\Omega_i)_{i \in V(G)}$ be a partition of $\Omega$ into parts of equal measure. We write $W_G$ for the graphon which is equal to $\ind{ij\in E(G)}$ on $\Omega_i\times\Omega_j$ and call it a \emph{graphon representation} of $G$.

A key insight of Garbe, Hladk\'{y}, \v{S}ileikis and Skerman is that the typical evolution of flip processes over quadratic timescales is highly concentrated along graphon-valued trajectories. This is formalized in their Transference Theorem, which we state below in a simplified form. It transfers problems about the typical evolution of random discrete-time graph-valued flip processes to problems about deterministic continuous-time graphon-valued trajectories. Consequently, understanding these graphon-valued trajectories is key to understanding the typical evolution of flip processes.

\begin{theorem}[Theorem 5.1 in~\cite{GarbeHladkySileikisSkerman}] \label{thm:transference}
For every rule $\rul$ there is a function $\trajrul{\rul}{}{} \colon \Gra \times [0,\infty) \to \Gra$ such that for every $\eps,T>0$ and flip process $(G_{\ell})_{\ell\in\N_0}$ with rule $\rul$ where $G_0$ has $n$ vertices, with probability at least $1-e^{-\Omega(n^2)}$ we have $\cutnd(W_{G_{\ell}},\trajrul{\rul}{}{}(W_{G_0},\frac{\ell}{n^2}))<\eps$ for all $0\le\ell\le Tn^2$.
\end{theorem}

We often write $\trajrul{\rul}{t}{W}$ instead of $\trajrul{\rul}{}{}(W,t)$ for the sake of brevity. We remark that in~\cite{GarbeHladkySileikisSkerman} the function $\trajrul{\rul}{}{}$ is constructed for each rule $\rul$ by first explicitly defining a \emph{velocity operator} $\vel[\rul] \colon \Gra \to \Kernel$ and then defining the \emph{trajectory}
$\left(\trajrul{\rul}{t}{W}\right)_{t\in[0,\infty)}$ starting at each $W\in\Gra$ to be the unique solution to the Banach-space-valued differential equation $\frac{\D}{\D t}\trajrul{\rul}{t}{W} = \vel[\rul]\trajrul{\rul}{t}{W}$ with the initial condition $\trajrul{\rul}{0}{W} = W$. More details are given in Section~\ref{sec:main-thm}.

The \emph{distribution} of the flip process started at a finite graph $G$ with rule $\rul$ is its usual probability distribution from the theory of stochastic processes. We say that two rules $\rul_1$ and $\rul_2$ of orders $k_1$ and $k_2$ respectively \emph{have the same flip process distributions} if for all finite graphs $G$ on $n\ge\max\{k_1,k_2\}$ vertices the flip process started at $G$ with rule $\rul_1$ has the same distribution as the flip process started at $G$ with rule $\rul_2$. Note that two rules with the same trajectories could have different flip process distributions. A rule $\rul$ is \emph{ignorant} if the replacement graph distribution $(\rul_{F,H})_{H\in\lgr{k}}$ is independent of the drawn graph $F\in\lgr{k}$. It was observed in~\cite[Section 4]{AraujoHladkyHngSileikis} that two ignorant rules $\rul_1$ and $\rul_2$ of the same order and with the same expected replacement edge count $d_{\rul} = \sum_{H\in\lgr{k}} \rul_{F,H}e(H)$ have the same trajectories, that is, $\trajrul{\rul_1}{}{} = \trajrul{\rul_2}{}{}$. More concretely, consider ignorant rules $\rul_1$ and $\rul_2$ of order $4$ where for $\rul_1$ the replacement graph is either complete or empty with probability $\frac{1}{2}$ each and for $\rul_2$ the replacement graph is a copy of $K_{1,3}$. Clearly, they have different replacement behaviour. On the other hand, both rules are ignorant rules of the same order with the same expected replacement edge count $d_{\rul_1} = d_{\rul_2} = 3$, so they have the same trajectories. Indeed, the lens of dense graph limits captures the macroscopic profile of flip processes but is insensitive to lower order fluctuations. This leads us to the following question.

\begin{question} \label{qn:rule-traj-equal}
When do two rules $\rul_1$ and $\rul_2$ have the same trajectories?
\end{question}

In this paper, we focus on the dense graph limits perspective of flip processes and fully resolve Question~\ref{qn:rule-traj-equal}, thereby providing a complete characterization of equivalence classes of flip process rules with the same trajectories.

We introduce some notation and terminology. An \emph{ordered set} is a set $S$ equipped with a linear order $\le_S$. For a set $S$ write $\binom{S}{2} = \{\{i,j\} \in S^2:i \ne j\}$ and $S^{(2)} = \{(i,j) \in S^2:i \ne j\}$. For the sake of brevity we often write $ij$ for an unordered pair $\{i,j\}$. A \emph{rooted graph} is a pair $\F = (F,R)$ consisting of a graph $F$ and an ordered set $R \subseteq V(F)$ of \emph{roots}. Set $V(\F) = V(F)$, $v(\F) = v(F)$, $E(\F) = E(F)$, $e(\F) = e(F)$, $R(\F) = R$, $r(\F) = |R|$, $V'(\F) = V(F) \sm R$ and $v'(\F) = v(F) - |R|$. When $R = (a,b)$ is an ordered pair, we often write $F^{a,b}$ for the rooted graph $(F,(a,b))$. For $k\in\N$ write $\cg_k := \lgr{k} \times [k]^{(2)}$ for the collection of rooted graphs $F^{a,b}$ where $F\in\lgr{k}$ and $(a,b)\in[k]^{(2)}$.

Let us explain the ideas motivating our main result. Let $\rul$ be a rule of order $k$. Heuristically, a graphon $W$ represents the adjacency matrix of a large graph $G$ on $n$ vertices and $W(x,y)$ represents the edge density of an auxiliary bipartite pair $(X,Y)$ of tiny subsets of $V(G)$ at $(x,y)\in\Omega^2$. Since $X$ and $Y$ are tiny, the rate of change $\vel[\rul]W(x,y)$ of the edge density of $(X,Y)$ is dominated by replacements where the drawn graph has exactly one vertex in each of $X$ and $Y$. We sample a copy of a labelled graph $F\in\lgr{k}$ with distinct vertices $a,b\in[k]$ in $X$ and $Y$ respectively with probability $|X||Y|n^{-2}\cdot\Troot{F^{a,b}}W(x,y)$. Here $\Troot{F^{a,b}}$ is an operator (see~\eqref{eq:operator-induced-density-rooted}) where $\Troot{F^{a,b}}W(x,y)$ is the induced density of $F$ in $W$ with the vertices $a$ and $b$ mapped to $x$ and $y$ respectively. Now the expected change in edge density on the pair $(X,Y)$ due to the replacement of $F$ is given by $\frac{q_{F^{a,b},\rul}}{|X||Y|}$ where $q_{F^{a,b},\rul} = \sum_{H\in\lgr{k}} \rul_{F,H}\cdot\ind{ab \in E(H)} - \ind{ab \in E(F)}$. Hence, by the law of total expectation, ignoring lower order terms, and rescaling time by $n^2$ because $(\Omega,\pi)$ is a probability space rather than a set of $n$ vertices, we obtain
\[\vel[\rul]W(x,y) = \sum_{F^{a,b}\in\cg_k} q_{F^{a,b},\rul}\cdot\Troot{F^{a,b}}W(x,y)\;.\]
Let us briefly comment on the $k=1$ case. One obviously cannot find distinct vertices $a,b\in[k]$, so the sum above is vacuous and trivially always zero; this is consistent with the observation that the only rule of order $1$ is a rule which does nothing.

Here $\Troot{F^{a,b}}W(x,y)$ is invariant under the relabelling of vertices of $F^{a,b}$ (including $a$ and $b$), which yields an isomorphic rooted graph in the sense of the definitions on page~\pageref{isomorphic}. Hence, we may focus on equivalence classes of rooted graphs $F^{a,b}$ related by relabelling; write $\cj_k$ for the set of such equivalence classes for $\cg_k$. For each $J\in\cj_k$ let $\Troot{J}W(x,y)$ be the unique value of $\Troot{F^{a,b}}W(x,y)$ for all $F^{a,b} \in J$ and set $q_{J,\rul} = \sum_{F^{a,b} \in J}q_{F^{a,b},\rul}$. Hence, we have 
\[\vel[\rul]W(x,y) = \sum_{J\in\cj_k}q_{J,\rul}\cdot\Troot{J}W(x,y)\;.\]

Now consider two rules $\rul_1$ and $\rul_2$ of the same order $k$. Suppose that $q_{J,\rul_1} = q_{J,\rul_2}$ for all $J\in\cj_k$. Then we have $\vel[\rul_1]W(x,y) = \vel[\rul_2]W(x,y)$ for all graphons $W$ and $(x,y)\in\Omega^2$, which implies that $\rul_1$ and $\rul_2$ have the same trajectories. Our main result is that the converse also holds, thereby giving a characterization of the equivalence classes of flip process rules of the same order with the same trajectories. More details are given in Section~\ref{sec:main-thm}.

\begin{theorem} \label{thm:rule-vel-char-intro}
The following are equivalent for rules $\rul_1$ and $\rul_2$ of the same order $k$.
\begin{enumerate}[label=(A\arabic{*})]
\item \label{item:same-coeff-intro} We have $q_{J,\rul_1} = q_{J,\rul_2}$ for all $J\in\cj_k$.
\item \label{item:same-traj-intro} $\trajrul{\rul_1}{}{} = \trajrul{\rul_2}{}{}$.
\end{enumerate}
\end{theorem}

Theorem~\ref{thm:rule-vel-char-intro} characterizes equivalence classes of rules of the same order with the same trajectories. Let us now demonstrate how it may be applied to determine whether two rules $\rul_1$ and $\rul_2$ of different orders have the same trajectories. Say the orders of $\rul_1$ and $\rul_2$ are $k_1$ and $k_2$ respectively; by symmetry we may assume $k_1<k_2$. Let $\rul^{*}_1$ be the rule of order $k_2$ where the replacement graph $H$ is obtained from the drawn graph $F$ by replacing the subgraph of $F$ induced on $[k_1]$ according to rule $\rul_1$ and leaving all other pairs unchanged. Since $\rul^{*}_1$ acts exclusively on pairs on $[k_1]$ and on those pairs it acts in accordance with $\rul_1$, it follows that $\rul^{*}_1$ and $\rul_1$ have the same trajectories. Let $F^{a,b}$ be a rooted graph with $F\in\lgr{k_2}$ and distinct $a,b\in[k_2]$. If $a,b\in[k_1]$, then we have $q_{F^{a,b},\rul^{*}_1} = q_{G^{a,b},\rul_1}$ with $G$ being the subgraph of $F$ induced on $[k_1]$. Otherwise, we have $q_{F^{a,b},\rul^{*}_1} = 0$. Now we apply Theorem~\ref{thm:rule-vel-char-intro} to determine whether $\trajrul{\rul^{*}_1}{}{} = \trajrul{\rul_2}{}{}$, which is equivalent to $\trajrul{\rul_1}{}{} = \trajrul{\rul_2}{}{}$. All in all, this gives us a complete characterization of equivalence classes of flip process rules with the same trajectories.

Many combinatorially interesting rules turn out to be symmetric deterministic rules. Write $S_k$ for the symmetric group on $[k]$. We have the following natural action of $S_k$ on $\lgr{k}$. Given $\sigma \in S_k$ and $F\in\lgr{k}$ we define $\sigma \cdot F$ to be the graph with vertex set $[k]$ and edge set $\{\sigma(i)\sigma(j):ij\in E(F)\}$. A rule $\rul$ of order $k$ is \emph{symmetric} if for all $F,H\in\lgr{k}$ and $\sigma \in S_k$ we have $\rul_{\sigma \cdot F,\sigma \cdot H} = \rul_{F,H}$ and \emph{deterministic} if for each $F\in\lgr{k}$ there is a unique $G\in\lgr{k}$ such that $\rul_{F,H} = \ind{H=G}$. Indeed, naturally arising rules are likely to act symmetrically on isomorphic graphs. For example, naturally arising rules of order $3$ are likely to treat the cherry graph $K_{1,2}$ symmetrically, no matter whether the middle vertex is labelled $1$, $2$ or $3$. The class of symmetric deterministic rules contains many natural families of rules which were highlighted in~\cite{GarbeHladkySileikisSkerman} and studied in~\cite{AraujoHladkyHngSileikis}. These include the complementing rules, the component completion rules, the extremist rules, and the clique removal rules.

Another interesting class of rules is the collection of \emph{unique} rules, that is, rules such that no other rule of the same order has the same trajectories. As an application of our main theorem, we prove that the unique rules are precisely the symmetric deterministic rules and the rules of order $2$. Observe that the rules of order $2$ form a very simple family of symmetric rules which is parametrized by just two parameters: the probabilities of replacing an edge by a non-edge and a non-edge by an edge.

\begin{theorem} \label{thm:rule-unique-symm-det}
A rule $\rul$ is unique if and only if it is a symmetric deterministic rule or a rule of order $2$.
\end{theorem}

The remainder of this paper is organized as follows. In Section~\ref{sec:main-thm} we formally state our main theorem, Theorem~\ref{thm:rule-vel-char}, after introducing the necessary notation, terminology and concepts. In Section~\ref{sec:prove-thm} we state the main lemmas for our proof of Theorem~\ref{thm:rule-vel-char} and give our proof of Theorem~\ref{thm:rule-vel-char}. We provide the proofs of the main lemmas in Section~\ref{subsec:proofs-lemmas}. In Section~\ref{sec:prove-thm-unique} we provide our proof of Theorem~\ref{thm:rule-unique-symm-det}.

\section{Main Theorem} \label{sec:main-thm}

In this section we give the full statement of our main result. We begin with some definitions we need to formally introduce the concept of a graphon-valued trajectory. For a rooted graph $F^{a,b}$ we define an operator $\Troot{F^{a,b}}\colon\Kernel\rightarrow L^\infty(\Omega^2)$ by setting
\begin{equation} \label{eq:operator-induced-density-rooted}
\Troot{F^{a,b}}W(x,y) \coloneqq \int \prod_{ij \in E(F)} W(x_i,x_j) \prod_{ij \notin E(F)} \left( 1 - W(x_i,x_j) \right) \D \pi^{V'(F^{a,b})}\;,
\end{equation}
with $x_a = x$ and $x_b = y$ fixed and the integral taken over $(x_i)_{i\in V'(F^{a,b})}\in\Omega^{V'(F^{a,b})}$. The \emph{velocity operator} for a rule $\rul$ of order $k$ is the operator $\vel[\rul]\colon\Kernel\to\Kernel$ given by
\begin{equation} \label{eq:velocity}
\vel[\rul]W := \sum_{F^{a,b}\in\cg_k} \left(\sum_{H\in\lgr{k}} \rul_{F,H} \cdot \ind{ab \in E(H)} - \ind{ab \in E(F)}\right) \Troot{F^{a,b}}W \;.
\end{equation}

Now we state the definition of a graphon-valued trajectory.

\begin{definition}[{Definition~4.4 in~\cite{GarbeHladkySileikisSkerman}}] \label{def:traj}
Given a rule $\rul$, a \emph{trajectory} starting at a kernel $W \in \Kernel$ is a differentiable function $\trajrul{\rul}{\cdot}{W} \colon I \to (\Kernel, \Linf{\cdot})$ defined on an open interval $I \subseteq \R$ containing zero which satisfies the autonomous differential equation
\[ \frac{\D}{\D t}\trajrul{\rul}{t}{W} = \vel[\rul]\trajrul{\rul}{t}{W} \]
with the initial condition $\trajrul{\rul}{0}{W} = W$.
\end{definition}

The following theorem guarantees that there is a unique trajectory with maximal domain starting at each kernel and establishes some useful facts about trajectories.

\begin{theorem}[{Theorem~4.5 in~\cite{GarbeHladkySileikisSkerman}}] \label{thm:traj}
The following hold for any rule $\rul$ and any kernel $W\in\Kernel$.
\begin{enumerate}[label=(\roman{*})]
\item \label{item:traj-exist} There is an open interval $\mdom{\rul,W} \subseteq \R$ containing $0$ and a trajectory $\trajrul{\rul}{\cdot}{W} \colon \mdom{\rul,W} \to \Kernel$ starting at $W$ such that any other trajectory starting at $W$ is a restriction of $\trajrul{\rul}{\cdot}{W}$ to a subinterval of $\mdom{\rul,W}$.
\item \label{item:traj-semi} For any $u\in\mdom{\rul,W}$ we have $\mdom{\rul,\traj{u}{W}} = \{ t \in \R : t + u \in \mdom{\rul,W}\}$ and for every $t\in\mdom{\rul,\traj{u}{W}}$ we have $\trajrul{\rul}{t}{\trajrul{\rul}{u}{W}} = \trajrul{\rul}{t+u}{W}$.
\item \label{item:traj-life}
If $W\in\Gra$ is a graphon, then the set $\life{\rul,W} := \{ t \in \mdom{\rul,W} : \trajrul{\rul}{t}{W} \in \Gra\}$ is a closed interval containing $[0,\infty)$.
\end{enumerate}
\end{theorem}

We introduce some further definitions. We say that rooted graphs $\F$ and $\HH$ are \emph{isomorphic}\label{isomorphic} if there is an \emph{isomorphism} from $\F$ to $\HH$, that is, a bijective function $\phi \colon V(\F) \to V(\HH)$ such that
\begin{enumerate}[label=(F\arabic{*})]
\item \label{item:fn-relpres} for all $u,v\in V(\F)$ we have $uv\in E(\F)$ if and only if $\phi(u)\phi(v)\in E(\HH)$,
\item \label{item:fn-rootresp} for all $v\in V(\F)$ we have $\phi(v) \in R(\HH)$ if and only if $v \in R(\F)$, and
\item \label{item:fn-rootordpres} for all $u,v \in R(\F)$ we have $\phi(u)\le_{R(\HH)}\phi(v)$ if and only if $u\le_{R(\F)}v$;
\end{enumerate}
write $\F\cong\HH$. Let $k\in\N$. Write $\cj_k$ for the set of isomorphism classes of the collection $\cg_k$. For $J\in\cj_k$ and a rule $\rul$ of order $k$, set
\begin{equation} \label{eq:rule-orbit-coeff}
q_{J,\rul} = \sum_{F^{a,b}\in J} \left( \sum_{H\in\lgr{k}} \rul_{F,H}\cdot\ind{ab \in E(H)} - \ind{ab \in E(F)}\right)\;.
\end{equation}
We say that two rules $\rul_1$ and $\rul_2$ \emph{have the same trajectories} if for all $W\in\Kernel$ we have $\mdom{\rul_1,W} = \mdom{\rul_2,W}$ and furthermore for all $t\in\mdom{\rul_1,W}$ we have $\trajrul{\rul_1}{t}{W} = \trajrul{\rul_2}{t}{W}$. Our main theorem, which we state below, characterizes rules of the same order with the same trajectories. Together with the `lifting' procedure relating rules of different orders described just below Theorem~\ref{thm:rule-vel-char-intro}, this fully resolves Question~\ref{qn:rule-traj-equal}.

\begin{theorem} \label{thm:rule-vel-char}
The following are equivalent for rules $\rul_1$ and $\rul_2$ of the same order $k\in\N$.
\begin{enumerate}[label=(S\arabic{*})]
\item \label{item:same-coeff} For all $J\in\cj_k$ we have $q_{J,\rul_1} = q_{J,\rul_2}$\;.
\item \label{item:same-vel-graphon} $\vel[\rul_1] = \vel[\rul_2]$ on $\Gra$\;.
\item \label{item:same-vel-kernel} $\vel[\rul_1] = \vel[\rul_2]$ on $\Kernel$\;.
\item \label{item:same-traj} $\rul_1$ and $\rul_2$ have the same trajectories.
\end{enumerate}
\end{theorem}

\section{Proving Theorem~\ref{thm:rule-vel-char}} \label{sec:prove-thm}

In this section we prove Theorem~\ref{thm:rule-vel-char}. We briefly describe our proof strategy here. Our proof scheme is \ref{item:same-coeff} $\Rightarrow$ \ref{item:same-vel-kernel} $\Rightarrow$ \ref{item:same-traj} $\Rightarrow$ \ref{item:same-vel-graphon} $\Rightarrow$ \ref{item:same-coeff}. We note that rooted induced densities are invariant under rooted graph isomorphism, so for $J\in\cj_k$ we may set $\Troot{J}$ to be the unique operator equal to $\Troot{F^{a,b}}$ for all $F^{a,b}\in J$ and the first implication follows. The second implication follows from the definition and uniqueness of trajectories (Definition~\ref{def:traj} and Theorem~\ref{thm:traj}), while the third implication follows from the smoothness of the velocity operator (Lemmas~\ref{lem:traj-space-perturb} and~\ref{lem:traj-time-perturb}). Indeed, the final implication \ref{item:same-vel-graphon} $\Rightarrow$ \ref{item:same-coeff} is the most difficult step. The goal here is to show that the collection $\ci = \{\Troot{J}:J\in\cj_k\}$ of operators is linearly independent.

We apply a couple of key insights to achieve this goal. Firstly, we focus on parametrized graphon representations of graphs because evaluating induced densities at these graphons nicely discretizes them into sums of monomials indexed over functions between rooted graphs that satisfy~\ref{item:fn-relpres}, \ref{item:fn-rootresp} and~\ref{item:fn-rootordpres}. Secondly, we focus on twinfree versions of rooted graphs because Lemma~\ref{lem:function-twins-tf-inj} tells us that the aforementioned functions are necessarily injective when the graphs involved are twinfree; in particular, this induces a partial order on the set of twinfree versions of rooted graphs. Combining these two insights leads us to Lemma~\ref{lem:rooted-density-quantify}, which provides an explicit quantification of induced densities at parametrized graphon representations of graphs with a unique monomial assigned to each relevant rooted graph. Finally, the desired linear independence follows from the linear independence of monomials in polynomial functions (see Lemma~\ref{lem:nonzero-polynomial-witness}).

The remainder of this section is organized into subsections as follows. First, we introduce the key concepts and associated key lemmas in Section~\ref{subsec:concepts-lemmas}. Next, we prove Theorem~\ref{thm:rule-vel-char} in Section~\ref{subsec:proof-thm}. Finally, we prove the key lemmas in Section~\ref{subsec:proofs-lemmas}.

\subsection{Key concepts and main lemmas} \label{subsec:concepts-lemmas}

In this subsection we introduce some notation, the key concepts and our main lemmas. The proofs of the lemmas are given later in Section~\ref{subsec:proofs-lemmas}. 

Let $\F$ be a rooted graph. Write $N_\F(v;S)$ for the set of neighbours of $v\in V(\F)$ in $\F$ which lie in $S\subseteq V(\F)$. We omit the set $S$ when we mean $S=V(\F)$, and we often omit $\F$ when it is clear from context. We say that $u,v\in V(\F)$ are \emph{twins} if $N(u) = N(v)$. Write $U(\F)$ for the set of vertices in $V'(\F)$ with a twin in $R(\F)$. Set $V^*(\F) = V'(\F) \sm U(\F)$, $u(\F) = |U(\F)|$ and $v^*(\F) = v'(\F) - u(\F)$. To give an example, obtain a rooted graph $\F$ by taking a copy of $K_4$ on $[4]$, deleting the edge $34$ and setting $(2,3)$ to be the ordered pair of roots. Here $V'(\F) = \{1,4\}$ is the set of non-roots. The only twins are $3$ and $4$, so $U(\F) = \{4\}$ and $V^*(\F)=\{1\}$. We say that $\F$ is \emph{twinfree} if for all $u,v\in V(\F)$ we have that $|R(\F)\cap\{u,v\}| = 1$ or that $N(u) = N(v)$ implies $u = v$. In particular, we permit pairs of twins with one root vertex and one non-root vertex. For a rooted graph $\F=(F,R)$ write $\F^\times$ for the graph $F$ and $\F^\nth$ for the rooted graph $(F,\nth)$.

Our first lemma deals with twin vertices in the context of functions between rooted graphs that satisfy~\ref{item:fn-relpres} and~\ref{item:fn-rootresp}.

\begin{lemma} \label{lem:function-twins-tf-inj}
The following hold for any function $\phi\colon V(\F)\to V(\HH)$ from a rooted graph $\F$ to a rooted graph $\HH$ that satisfies~\ref{item:fn-relpres} and~\ref{item:fn-rootresp}.
\begin{enumerate}[label=(\roman{*})]
\item \label{item:function-twins} If $u,v\in V(\F)$ are such that $\phi(u)$ and $\phi(v)$ are twins in $\HH$, then $u$ and $v$ are twins in $\F$.
\item \label{item:function-tf-inj} If $\F$ is twinfree, then $\phi$ is injective.
\end{enumerate}
\end{lemma}

Next, we describe a pair of complementary operations on rooted graphs. The first one shrinks clusters of twins in a rooted graph to obtain a twinfree version, while the second one blows up vertices of a (twinfree) rooted graph to clusters of twins. Let $(S,\le_S)$ be a finite ordered set with a partition $\cs = \{S_i\}_{i\in I}$ indexed by a finite set $I$. Writing $a_i$ for the least element of $S_i$ according to $\le_S$, the linear order $\le_I$ on $I$ \emph{induced by $\le_S$} is given by $i\le_Ij$ if and only if $a_i\le_Sa_j$. For a rooted graph $\F$ we define an equivalence relation $\sim_{\tf(\F)}$ on $V(\F)$ as follows. For $u,v\in V(\F)$ we have $u\sim_{\tf(\F)}v$ if $u=v$ or both $N(u) = N(v)$ and $|R(\F)\cap\{u,v\}| \ne 1$ hold. Write $\cc(\F)$ for the set of equivalence classes of $\sim_{\tf(\F)}$. For $v\in V(\F)$ we write $v_{\tf(\F)}$ for the unique equivalence class of $\sim_{\tf(\F)}$ containing $v$. For $A\subseteq V(\F)$ we set $A_{\tf(\F)} = \{a_{\tf(\F)}:a\in A\}$. The \emph{twinfree version} of $\F$ is the rooted graph $\tf(\F)$ on $\cc(\F)$ whose edges are the pairs $IJ$ for which there are $i\in I$ and $j\in J$ such that $ij\in E(\F)$ and whose roots are the elements of $R(\F)_{\tf(\F)}$ with the linear order induced by $\le_{R(\F)}$. That $\tf(\F)$ is twinfree follows easily from the definitions.

For a rooted graph $\F$ and a vector $\m = (m_i)_{i \in V(\F)} \in \N_0^{V(\F)}$ the \emph{$\m$-blowup} of $\F$ is the rooted graph $\F(\m)$ obtained as follows. First, replace each $i\in V(\F)$ by an independent set $S_i$ of $m_i$ vertices and each edge $ij\in E(\F)$ by a complete bipartite graph on $(S_i,S_j)$. We say that we \emph{blow up} each vertex $i \in V(\F)$ by a factor of $m_i$ to obtain $\F(\m)$ from $\F$ and often omit mention of vertices $i\in V(\F)$ where $m_i = 1$. Next, let the set of roots be $\bigcup_{i\in R(\F)}S_i$ ordered so that each $S_i$ forms an interval and for all $u\in S_i$ and $v\in S_j$ with $i\le_{R(\F)}j$ we have $u\le_{R(\F(\m))}v$. Note that the linear order is not specified within each interval $S_i$; any order on each interval $S_i$ would give the same rooted graph up to isomorphism, so we conveniently avoid making these arbitrary choices. It follows easily from the definitions that every rooted graph $\F$ with $r(\F)\le2$ is isomorphic to the $\m$-blowup of its twinfree version $\tf(\F)$ for some vector $\m$ of positive integers. Since in practice our rooted graphs have at most two roots, this is sufficient for our purposes; we remark that for many rooted graphs with at least three roots the statement above actually does not hold due to difficulties in handling the order on the roots.

Now we define notation and terminology to describe our quantification of induced densities at certain graphons. Denote by $Z(G)$ the collection of vectors $\z = (z_i)_{i \in V(G)} \in [0,1]^{V(G)}$ satisfying $\sum_{i \in V(G)}z_i = 1$. For $\z\in Z(G)$ write $\Om(\z)$ for a partition $(\Omega_i)_{i \in V(G)}$ of $\Omega$ such that $\Omega_i$ has measure $z_i$. For $\z\in Z(G)$ and $\Om(\z)$ we say that the graphon $W_G^{\z}$ which is equal to $\ind{ij\in E(G)}$ on $\Omega_i\times\Omega_j$ is a \emph{$\z$-scaled graphon representation of $G$}. For a graph $G$ and an ordered subset $X \subseteq V(G)$ the \emph{$X$-rooted version of $G$} is the rooted graph $\F$ obtained from $G$ by blowing up each vertex in $X$ by a factor of $2$ and designating as the ordered set of roots the set of duplicate vertices with the linear order induced by $\le_X$. For $(x,y)\in\Omega^2$ the \emph{$(x,y)$-rooted version of $G$ for partition $(\Omega_i)_{i \in V(G)}$} is the $X$-rooted version of $G$ where the ordered subset $X=\{u,v\}\subseteq V(G)$ is such that $x\in\Omega_u$, $y\in\Omega_v$ and $u\le_X v$. Note that we may have $u=v$; in this case, $X$ is a singleton set. Finally, we write $x_{\rt}$ (resp.\ $y_{\rt}$) for the root vertex whose twin $w$ in $X$ satisfies $x\in\Omega_w$ (resp.\ $y\in\Omega_w$).

The following lemma gives an explicit quantification of induced densities at parametrized graphon representations of graphs. The availability of such an explicit quantification is central to our proof of Theorem~\ref{thm:rule-vel-char}. A function $\phi \colon V(\F) \to V(\HH)$ from a rooted graph $\F$ to a rooted graph $\HH$ is \emph{induced copy} if it satisfies~\ref{item:fn-relpres}, \ref{item:fn-rootresp} and~\ref{item:fn-rootordpres}.

\begin{lemma} \label{lem:rooted-density-quantify}
Let $\F$ and $\G$ be twinfree rooted graphs with $R(\G)=\nth$. Let $\z\in Z(\G^\times)$ and $W_{\G^\times}^{\z}$ be a $\z$-scaled graphon representation of the graph $\G^\times$ with partition $\Om = \Om(\z)$. Let $(x,y)\in\Omega^2$ and let $\HH$ be the $(x,y)$-rooted version of $\G^\times$ for partition $\Om$. The following hold for all vectors $\m = (m_i)_{i \in V(\F)} \in \N^{V(\F)}$ satisfying $\sum_{i\in R(\F)}m_i = 2$.
\begin{enumerate}[label=(\roman{*})]
    \item \label{item:rooted-density-strict} If $r(\F) > r(\HH)$ or both $r(\F) = r(\HH)$ and $v^*(\F) > v^*(\HH)$ hold, then $\Troot{\F(\m)}W_{\G^\times}^{\z}(x,y) = 0.$
    \item \label{item:rooted-density-equal} If $r(\F) = r(\HH)$ and $v^*(\F) = v^*(\HH)$, then we have
    \begin{equation*}
        \Troot{\F(\m)}W_{\G^\times}^{\z}(x,y) = \sum_{\phi \colon V(\F) \to V(\HH)} \prod_{i \in V'(\F)} z_{\phi(i)}^{m_i}\;,
    \end{equation*}
    where the sum is over all induced copy functions $\phi \colon V(\F) \to V(\HH)$.
\end{enumerate}
\end{lemma}

The following lemmas give bounds on trajectories and the velocity operator. For $k\in\N$ set $C_k = 2^{\binom{k}{2}-1}k^2(k-1)^2$.

\begin{lemma}[{Lemma~4.8 in~\cite{GarbeHladkySileikisSkerman}}] \label{lem:traj-space-perturb}
Given a rule $\rul$ of order $k\in\N$ and $U,W\in\Gra$ we have
\[ \Linf{\vel[\rul] U - \vel[\rul] W } \le C_k \Linf{U - W}\;. \]
\end{lemma}

\begin{lemma}[{Lemma~4.15 in~\cite{GarbeHladkySileikisSkerman}}] \label{lem:traj-time-perturb}
Given a rule $\rul$ of order $k\in\N$, $W\in\Gra$ and $\delta > 0$ we have
\[ \Linf{\delta^{-1}(\traj{\delta}{W} - W) - \vel[\rul]W} \le C_kk(k-1) \delta\;. \]
\end{lemma}

Finally, we state a technical lemma which establishes the linear independence of monomials in polynomial functions. We first introduce the necessary notation, terminology and concepts. A $k$-variable \emph{polynomial function} is a function $f \colon \R^k \to \R$ given by
\begin{equation} \label{eq:poly-fn}
f(x_1,\dots,x_k) = \sum_{\m = (m_i)_{i\in[k]}\in\N_0^k} a_{\m}\prod_{i\in[k]}x_i^{m_i},
\end{equation}
where only finitely many of the coefficients $a_{\m}\in\R$ are nonzero. A $k$-variable polynomial function is \emph{nonzero} if it has at least one nonzero coefficient $a_{\m}$ and \emph{zero} otherwise. The \emph{degree} of a nonzero $k$-variable polynomial function $f$ is
\[\deg(f) := \max\left\{\sum_{i\in[k]}m_i:\m = (m_i)_{i\in[k]}\in\N_0^k,a_{\m}\ne0\right\}.\]
A \emph{zero} of a $k$-variable polynomial function $f$ is a vector $(x_i)_{i\in[k]}\in\R^k$ such that $f(x_1,\dots,x_k)=0$.

The following lemma, which is a special case of the Schwartz--Zippel Lemma~\cite{Schwartz,Zippel}, asserts that each nonzero multivariate polynomial function is not identically zero when taking values in sufficiently large sets.

\begin{lemma} \label{lem:nonzero-polynomial-witness}
For any $k\in\N$, finite subset $S\subseteq \R$ and nonzero $k$-variable polynomial function $f$ with $\deg(f)<|S|$, there exist $x_1,\dots,x_k \in S$ such that $f(x_1,\dots,x_k) \ne 0$.
\end{lemma}

\subsection{Proof of Theorem~\ref{thm:rule-vel-char}} \label{subsec:proof-thm}

In this subsection we prove Theorem~\ref{thm:rule-vel-char}. As mentioned at the beginning of this section, we use the proof scheme \ref{item:same-coeff} $\Rightarrow$ \ref{item:same-vel-kernel} $\Rightarrow$ \ref{item:same-traj} $\Rightarrow$ \ref{item:same-vel-graphon} $\Rightarrow$ \ref{item:same-coeff}. We shall present the proof of each implication separately.

\begin{proof}[Proof of \ref{item:same-coeff} $\Rightarrow$ \ref{item:same-vel-kernel}]
Let $F^{a,b},G^{c,d}\in J\in\cj_k$, $W\in\Kernel$ and $(x,y)\in\Omega^2$. Fix an isomorphism $\phi$ from $F^{a,b}$ to $G^{c,d}$. Set $x_a = x'_c = x$ and $x_b = x'_d = y$. By~\eqref{eq:operator-induced-density-rooted} we have
\begin{align*}
\Troot{F^{a,b}}W(x,y) & = \int \prod_{ij \in E(F)} W(x_i,x_j) \prod_{ij \notin E(F)} \left( 1 - W(x_i,x_j) \right) \D \pi^{V'(F^{a,b})} \\
& = \int \prod_{\phi(i)\phi(j) \in E(G)} W(x_i,x_j) \prod_{\phi(i)\phi(j) \notin E(G)} \left( 1 - W(x_i,x_j) \right) \D \pi^{V'(F^{a,b})} \\
& = \int \prod_{ij \in E(G)} W(x'_i,x'_j) \prod_{ij \notin E(G)} \left( 1 - W(x'_i,x'_j) \right) \D \pi^{V'(G^{c,d})} = \Troot{G^{c,d}}W(x,y)\;.
\end{align*}
Hence, for $J\in\cj_k$ we may write $\Troot{J}$ for the unique operator equal to $\Troot{\F}$ for all $\F\in J$. For a rule $\rul$ of order $k$ we rewrite~\eqref{eq:velocity} to obtain
\begin{equation} \label{eq:velocity-iso}
\vel[\rul]W = \sum_{\F\in\cg_k} q_{\F,\rul}\cdot\Troot{\F}W = \sum_{J\in\cj_k} q_{J,\rul}\cdot\Troot{J}W\;.
\end{equation}
It follows that~\ref{item:same-coeff} implies~\ref{item:same-vel-kernel}.
\end{proof}

\begin{proof}[Proof of \ref{item:same-vel-kernel} $\Rightarrow$ \ref{item:same-traj}]
Suppose that $\vel[\rul_1] = \vel[\rul_2]$ on $\Kernel$ and fix $W\in\Kernel$. By Theorem~\ref{thm:traj}\ref{item:traj-exist} we have trajectories $\trajrul{\rul_1}{\cdot}{W}$ and $\trajrul{\rul_2}{\cdot}{W}$ starting at $W$ which have maximal domains $\mdom{\rul_1,W}$ and $\mdom{\rul_2,W}$, respectively, that are open intervals containing $0$. The combination of~\ref{item:same-vel-kernel} and Definition~\ref{def:traj} implies that $\trajrul{\rul_1}{\cdot}{W}$ is a trajectory starting at $W$ for $\rul_2$ and $\trajrul{\rul_2}{\cdot}{W}$ is a trajectory starting at $W$ for $\rul_1$, so by Theorem~\ref{thm:traj}\ref{item:traj-exist} $\trajrul{\rul_1}{\cdot}{W}$ and $\trajrul{\rul_2}{\cdot}{W}$ are restrictions of each other. Hence, we have $\mdom{\rul_1,W} = \mdom{\rul_2,W}$ and $\trajrul{\rul_1}{t}{W} = \trajrul{\rul_2}{t}{W}$ for all $t\in\mdom{\rul_1,W}$, which means that $\rul_1$ and $\rul_2$ have the same trajectories.
\end{proof}

\begin{proof}[Proof of \ref{item:same-traj} $\Rightarrow$ \ref{item:same-vel-graphon}]
Suppose that there is $W\in\Gra$ such that $\vel[\rul_1]W \ne \vel[\rul_2]W$. Setting $\eps = \Linf{\vel[\rul_1]W - \vel[\rul_2]W} > 0$, take $0< \eta \le \frac{\eps}{3C_kk(k-1)}$. By the triangle inequality for $\Linf{\cdot}$ and Lemmas~\ref{lem:traj-space-perturb} and~\ref{lem:traj-time-perturb}, we have
\begin{align*}
\eps = \Linf{\vel[\rul_1]W - \vel[\rul_2]W} & \le \Linf{\eta^{-1}(\trajrul{\rul_1}{\eta}{W} - \trajrul{\rul_2}{\eta}{W})} + \sum_{i=1}^2 \Linf{\eta^{-1}(\trajrul{\rul_i}{\eta}{W} - W) - \vel[\rul_i]W} \\
& \le \eta^{-1}\Linf{\trajrul{\rul_1}{\eta}{W} - \trajrul{\rul_2}{\eta}{W}} + 2C_kk(k-1)\eta.
\end{align*}
Rearranging yields $\Linf{\trajrul{\rul_1}{\eta}{W} - \trajrul{\rul_2}{\eta}{W}} \ge \eta(\eps - 2C_kk(k-1)\eta) \ge \frac{\eta\eps}{3} > 0$, which implies $\trajrul{\rul_1}{\eta}{W} \ne \trajrul{\rul_2}{\eta}{W}$. Hence, $\rul_1$ and $\rul_2$ do not have the same trajectories.
\end{proof}

\begin{proof}[Proof of \ref{item:same-vel-graphon} $\Rightarrow$ \ref{item:same-coeff}]
For $\F\in\cg_k$ set $\alpha_\F = q_{\F,\rul_1} -  q_{\F,\rul_2}$. Suppose that~\ref{item:same-coeff} does not hold; fix $\G\in\cg_k$ satisfying $\alpha_\G\ne0$ which is minimal first in $r(\tf(\G))$ and then in $v^*(\tf(\G))$. Equip $R(\G)_{\tf(\G^\nth)}$ with the linear order induced by $\le_{R(\G)}$ and write $\HH$ for the $R(\G)_{\tf(\G^\nth)}$-rooted version of $\tf(\G^\nth)^\times$. That $a_{\tf(\G)} \mapsto a_{\tf(\G^\nth)}$ gives well-defined bijections from $R(\G)_{\tf(\G)}$ to $R(\G)_{\tf(\G^\nth)}$ and from $V^*(\G)_{\tf(\G)}$ to $V^*(\G)_{\tf(\G^\nth)}$ follows from the definitions, so $r(\HH) = r(\tf(\G))$ and $v^*(\HH) = v^*(\tf(\G))$.

We shall prove that~\ref{item:same-vel-graphon} does not hold. Let $\z\in Z(\tf(\G^\nth)^\times)$ and $\Om = \Om(\z) = (\Omega_i)_{i \in V(\tf(\G^\nth)^\times)}$. Take $(x,y)\in\Omega^2$ so that $R(\G)_{\tf(\G^\nth)} = \{v\in V(\tf(\G^\nth)):\{x,y\}\cap\Omega_v\ne\nth\}$ and the first element $w$ in the linear order induced by $\le_{R(\G)}$ satisfies $x\in\Omega_w$. In particular, $\HH$ is the $(x,y)$-rooted version of $\tf(\G^\nth)^\times$ for the partition $\Om$. By~\eqref{eq:velocity-iso} we have 
\begin{equation} \label{eq:velocity-diff}
\vel[\rul_1] W_{\tf(\G^\nth)^\times}^{\z}(x,y) - \vel[\rul_2] W_{\tf(\G^\nth)^\times}^{\z}(x,y) = \sum_{\F\in\cg_k} \alpha_{\F} \cdot \Troot{\F}W_{\tf(\G^\nth)^\times}^{\z}(x,y) \;.
\end{equation}

We handle the summands in~\eqref{eq:velocity-diff} with the following categorization of $\F$.
\begin{enumerate}[label=(C\arabic{*})]
\item \label{case:smaller-rooted} $r(\tf(\F)) < r(\HH)$, or $r(\tf(\F)) = r(\HH)$ and $v^*(\tf(\F)) < v^*(\HH)$ hold.
\item \label{case:equal-all} $r(\tf(\F)) = r(\HH)$ and $v^*(\tf(\F)) = v^*(\HH)$ hold.
\item \label{case:larger-rooted} $r(\tf(\F)) > r(\HH)$, or $r(\tf(\F)) = r(\HH)$ and $v^*(\tf(\F)) > v^*(\HH)$ hold.
\end{enumerate}
Summands of category~\ref{case:smaller-rooted} are trivial because $\alpha_{\F} = 0$ by choice of $\G$. For summands of categories~\ref{case:equal-all} and~\ref{case:larger-rooted} we focus on $\Troot{\F}W_{\tf(\G^\nth)^\times}^{\z}(x,y)$. Since $\tf(\F)$ and $\tf(\G^\nth)$ are twinfree, by Lemma~\ref{lem:rooted-density-quantify} we have $\Troot{\F}W_{\tf(\G^\nth)^\times}^{\z}(x,y) > 0$ only if $\F$ is in category~\ref{case:equal-all} and there is an induced copy function $\phi \colon V(\tf(\F)) \to V(\HH)$.

Let $\phi \colon V(\tf(\F)) \to V(\HH)$ be an induced copy function for some $\F$ in category~\ref{case:equal-all}. By adding dummy non-root twins of root vertices to $\tf(\F)$, we may extend $\phi$ (uniquely, by Lemma~\ref{lem:function-twins-tf-inj}\ref{item:function-twins}) to an isomorphism $\chi \colon V(\HH) \to V(\HH)$; equivalently, $\tf(\F)$ is a copy of $\HH$ with some non-root twins of root vertices deleted. It follows that $\F$ is an $\m$-blowup of $\HH$ for $\m\in\ca(\HH)$, where $\ca(\HH)$ is the set of vectors $\m=(m_i)_{i\in V(\HH)}\in\N_0^{V(\HH)}$ where $m_i$ is positive for $i\in V^*(\HH)\cup R(\HH)$ and we have $\sum_{i\in V'(\HH)}m_i = k-2$ and $\sum_{i\in R(\HH)}m_i = 2$. Hence, by~\eqref{eq:velocity-diff} and Lemma~\ref{lem:rooted-density-quantify} we obtain
\begin{equation} \label{eq:velocity-diff-poly}
    \vel[\rul_1] W_{\tf(\G^\nth)^\times}^{\z}(x,y) - \vel[\rul_2] W_{\tf(\G^\nth)^\times}^{\z}(x,y) = \sum_{\m\in\ca(\HH)} \prod_{i\in V'(\HH)} z_i^{m_i} \cdot \alpha_{\HH(\m)}\;.
\end{equation}

Fix $q\in V'(\HH)$ and set $S=V'(\HH)$, $s=|S|$, $S'=S\sm\{q\}$ and $s'=s-1$. Since $\z$ is subject to the condition $\sum_{i\in S}z_i = 1$, we replace $z_q$ with $1-\sum_{i\in S'}z_i$ in~\eqref{eq:velocity-diff-poly} to obtain an $s'$-variable polynomial function $p$. The coefficient of each monomial $\prod_{i \in S'} z_i^{m_i}$ in $p$ has the form $\sum_{\n} C_{\m,\n} \cdot \alpha_{\HH(\n)}$ where the constants $C_{\m,\n}$ satisfy $C_{\m,\m} \ne 0$ and the sum is over all vectors $\n$ such that $n_i \le m_i$ for all $i\in S'$. Since $\alpha_{\G} \ne 0$, we may pick $\p$ such that $\alpha_{\HH(\p)} \ne 0$ and $\alpha_{\HH(\n)} = 0$ for all $\n \ne \p$ where we have $n_i \le p_i$ for all $i\in V(\HH)\sm\{q\}$. For example, we may start with $\p$ such that $\G=\HH(\p)$, replace the current $\p$ with any $\n\ne\p$ such that $\alpha_{\HH(\n)}\ne0$ and $n_i \le p_i$ for all $i\in V(\HH)\sm\{q\}$ if possible, and take the vector $\p$ at termination. The coefficient of $\prod_{i \in S'} z_i^{p_i}$ in $p$ is the nonzero $C_{\p,\p} \cdot \alpha_{\HH(\p)}$, so $p$ is nonzero. Now $p$ is not identically zero on $[0,s^{-1}]^{s-1}$ by Lemma~\ref{lem:nonzero-polynomial-witness}, so we may pick $\y\in Z(\tf(\G^\nth)^\times)$ such that $\vel[\rul_1] W_{\tf(\G^\nth)^\times}^{\y} \ne \vel[\rul_2] W_{\tf(\G^\nth)^\times}^{\y}$.
\end{proof}

\subsection{Proofs of lemmas} \label{subsec:proofs-lemmas}

In this subsection we give the proofs of our lemmas from Section~\ref{subsec:concepts-lemmas}.

We introduce some further definitions. Let $\F$ and $\HH$ be rooted graphs. A function $\phi \colon V(\F) \to V(\HH)$ is \emph{relation-preserving} if it satisfies~\ref{item:fn-relpres}, \emph{root-respecting} if it satisfies~\ref{item:fn-rootresp}, and \emph{root-order-preserving} if it satisfies~\ref{item:fn-rootordpres}. In particular, an induced copy function between rooted graphs is relation-preserving, root-respecting and root-order-preserving.

\begin{proof}[Proof of Lemma~\ref{lem:function-twins-tf-inj}]
For part~\ref{item:function-twins} take $u,v\in V(\F)$ such that $\phi(u)$ and $\phi(v)$ are twins in $\HH$; by definition we have $N_\HH(\phi(u)) = N_\HH(\phi(v))$. Since $\phi$ is relation-preserving, we have $N_\F(u) = N_\F(v)$, that is, $u$ and $v$ are twins in $\F$. For part~\ref{item:function-tf-inj} take $u,v\in V(\F)$ satisfying $\phi(u) = \phi(v)$. By part~\ref{item:function-twins} $u$ and $v$ are twins in $\F$. Since $\F$ is twinfree and we cannot have $|R(\F)\cap\{u,v\}| = 1$ because $\phi$ is root-respecting, we must have $u = v$.
\end{proof}

\begin{proof}[Proof of Lemma~\ref{lem:rooted-density-quantify}]
For $\x = (x_i)_{i\in V(\F(\m))} \in \Omega^{V(\F(\m))}$ set 
\begin{equation*}
a(\x) = \prod_{ij \in E(\F(\m))} W_{\G^\times}^{\z}(x_i,x_j) \prod_{ij \notin E(\F(\m))} \left( 1 - W_{\G^\times}^{\z}(x_i,x_j) \right).
\end{equation*}
For $\x\in\Omega^{V(\F(\m))}$ define $f_{\x} \colon V(\F(\m)) \to V(\G^\times)$ as follows. For $i\in V(\F(\m))$ set $f_{\x}(i) = w$ where $x_i \in \Omega_w$. Let $(p,q) = R(\F(\m))$. Set $\cp = \{\x\in\Omega^{V(\F(\m))}:f_{\x}\textrm{ is relation-preserving}\}$ and $\cq = \{\x\in\cp:x_p=x,x_q=y\}$. For $\x\in\cq$ define functions $g_{\x} \colon V(\F(\m)) \to V(\HH)$ and $h_{\x} \colon V(\F) \to V(\HH)$ as follows. For $i\in V(\F(\m))$ set
\begin{equation*}
g_{\x}(i) = 
\begin{cases}
f_{\x}(i) & \textrm{ if }i \ne p,q \\
x_{\rt} & \textrm{ if }i = p \\
y_{\rt} & \textrm{ if }i = q.
\end{cases}
\end{equation*}
For $i\in V(\F)$ pick a representative $a_i \in S_i$ and set $h_{\x}(i) = g_{\x}(a_i)$. The following claim establishes some useful properties. 

\begin{claim} \label{claim:function-properties}
The following hold.
\begin{enumerate}[label=(\roman{*})]
\item \label{item:characterize-integrand-p} For all $\x\in\Omega^{V(\F(\m))}$ we have $a(\x) = \ind{\x\in\cp}$.
\item \label{item:rooted-projection-properties} For all $\x\in\cq$ the functions $g_{\x}$ and $h_{\x}$ are induced copy.
\item \label{item:strict-q-empty} If $\cq\neq\nth$ then either $r(\F) < r(\HH)$ or both $r(\F) = r(\HH)$ and $v^*(\F) \le v^*(\HH)$ hold.
\end{enumerate}
The following hold when $r(\F) = r(\HH)$ and $v^*(\F) = v^*(\HH)$.
\begin{enumerate}[label=(\roman{*}),resume]
\item \label{item:projection-blowup-equal-unique} For all $\x\in\cq$, $\ell\in V(\F)$ and $i,j\in S_{\ell}$ we have $g_{\x}(i) = g_{\x}(j)$. In particular, the function $h_{\x}$ is uniquely defined.
\item \label{item:projection-characterize} For every induced copy function $h \colon V(\F) \to V(\HH)$ and vector $\x\in\Omega^{V(\F(\m))}$, the following are equivalent.
\begin{enumerate}[label=(\alph{*})]
\item \label{subitem:proj-char-target} $\x\in\cq$ and $h_{\x} = h$.
\item \label{subitem:proj-char-characterization} $x_p = x$, $x_q = y$, and for all $\ell\in V'(\F)$ and $i\in S_{\ell}$ we have $x_i \in \Omega_{h(\ell)}$. 
\end{enumerate}
\end{enumerate}
\end{claim}

\begin{claimproof}
For part~\ref{item:characterize-integrand-p} note that the definitions of $W_{\G^\times}^{\z}$ and $f_{\x}$ imply that for all $\x\in\Omega^{V(\F(\m))}$ and all $i,j\in V(\F(\m))$ we have $W_{\G^\times}^{\z}(x_i,x_j) = \ind{f_{\x}(i)f_{\x}(j)\in E(\G^\times)}$. Hence, we have
\[ a(\x) = \prod_{i,j\in V(\F(\m))} \ind{ij\in E(\F(\m)) \Leftrightarrow f_{\x}(i)f_{\x}(j)\in E(\G^\times)} = \ind{\x\in\cp}. \]

For part~\ref{item:rooted-projection-properties} note that $g_{\x}$ is root-respecting and root-order-preserving by definition. Since $f_{\x}$ is relation-preserving for $\x\in\cq\subseteq\cp$ and $x_{\rt}$ and $y_{\rt}$ are copies of $f_{\x}(p)$ and $f_{\x}(q)$ respectively, it follows that $g_{\x}$ is relation-preserving. The function $h_{\x}$ inherits the desired properties from $g_{\x}$.

For part~\ref{item:strict-q-empty} suppose that there is $\x\in\cq$. Since $\F$ is twinfree, by Lemma~\ref{lem:function-twins-tf-inj}\ref{item:function-tf-inj} the function $h_{\x}$ is injective. Hence, we have $r(\F) \le r(\HH)$. Focusing on the case $r(\F) = r(\HH)$, by Lemma~\ref{lem:function-twins-tf-inj}\ref{item:function-twins} we also have $\phi(v) \in V^*(\HH)$ for all $v\in V^*(\F)$, so the injectivity of $h_{\x}$ gives $v^*(\F) \le v^*(\HH)$.

For part~\ref{item:projection-blowup-equal-unique} fix $\x\in\cq$ and take $\ell\in V(\F)$ and $i,j\in S_{\ell}$. By the definition of $\F(\m)$ we have $N_{\F(\m)}(i) = N_{\F(\m)}(j)$ and by part~\ref{item:rooted-projection-properties} $g_{\x}$ is relation-preserving and root-respecting, so we have $N_{\HH}(g_{\x}(i);\im(g_{\x})) = N_{\HH}(g_{\x}(j);\im(g_{\x}))$. Since $\F$ is twinfree, by Lemma~\ref{lem:function-twins-tf-inj}\ref{item:function-tf-inj} the function $h_{\x}$ is injective, which implies $|h_{\x}(R(\F))| \ge r(\F) \ge r(\HH)$ and $|h_{\x}(V^*(\F))| \ge v^*(\F) \ge v^*(\HH)$. But $h_{\x}$ is root-respecting, so we have $h_{\x}(R(\F)) = R(\HH)$. Furthermore, by Lemma~\ref{lem:function-twins-tf-inj}\ref{item:function-twins} applied to $h_{\x}$ we have $h_{\x}(V^*(\F)) = V^*(\HH)$. Now we have $h_{\x}(V^*(\F) \cup R(\F)) = V^*(\HH) \cup R(\HH)$ and $\im(h_{\x}) \subseteq \im(g_{\x})$, so every vertex in $V(\HH)\sm\im(g_{\x})$ has a twin in $\im(g_{\x})$. Hence, we have $N_{\HH}(g_{\x}(i)) = N_{\HH}(g_{\x}(j))$. But $\HH$ is twinfree, so we have $g_{\x}(i) = g_{\x}(j)$. In particular, the function $h_{\x}$ obtained is independent of the specific choices of the representatives $a_i \in S_i$.

Finally, we consider part~\ref{item:projection-characterize}. Let $h \colon V(\F) \to V(\HH)$ be an induced copy function. Observe that the definitions of $\cq$ and $h_{\x}$ imply that any $\x$ satisfying the conditions in~\ref{subitem:proj-char-target} also satisfies the conditions in~\ref{subitem:proj-char-characterization}. Now suppose that $\x$ satisfies the conditions in~\ref{subitem:proj-char-characterization}. Let $i \in V(\F(\m))$ and take the unique $\ell\in V(\F)$ so that $i\in S_\ell$; by assumption we have $f_{\x}(i) = h(\ell)$. Since $h$ is relation-preserving, the definitions of $\HH$ and of blowups imply that $f_{\x}$ is relation-preserving, so we have $\x\in\cq$. Now the definition of $h_{\x}$ and its uniqueness by~\ref{item:projection-blowup-equal-unique} imply that $h_{\x}=h$, completing the proof that $\x$ satisfies the conditions in~\ref{subitem:proj-char-target}.
\end{claimproof}

We now apply Claim~\ref{claim:function-properties} to prove parts~\ref{item:rooted-density-strict} and~\ref{item:rooted-density-equal}. By~\eqref{eq:operator-induced-density-rooted} and Claim~\ref{claim:function-properties}\ref{item:characterize-integrand-p} we have
\begin{equation*}
\Troot{\F(\m)}W_{\G^\times}^{\z}(x,y) = \int_{\y \in \Omega^{V'(\F(\m))}} \ind{(\y,x,y)\in\cq} \D \pi^{V'(\F(\m))}\;.
\end{equation*}
Note that Claim~\ref{claim:function-properties}\ref{item:strict-q-empty} now immediately implies part~\ref{item:rooted-density-strict}. For part~\ref{item:rooted-density-equal}, observe that since $h_{\x}$ is uniquely defined for $x\in\cq$ by Claim~\ref{claim:function-properties}\ref{item:projection-blowup-equal-unique} and Claim~\ref{claim:function-properties}\ref{item:projection-characterize} gives an exact characterization of the partition of $\cq$ according to the associated induced copy function $h_{\x}$, we obtain
\begin{equation*}
\Troot{\F(\m)}W_{\G^\times}^{\z}(x,y) = \int_{\y \in \Omega^{V'(\F(\m))}} \ind{(\y,x,y)\in\cq} \D \pi^{V'(\F(\m))} = \sum_{\phi} \prod_{i \in V'(\F)} z_{\phi(i)}^{m_i}\;,
\end{equation*}
where the sum is over all induced copy functions $\phi \colon V(\F) \to V(\HH)$.
\end{proof}

\section{Proving Theorem~\ref{thm:rule-unique-symm-det}} \label{sec:prove-thm-unique}

In this section we prove Theorem~\ref{thm:rule-unique-symm-det}. We briefly describe our proof strategy here. Theorem~\ref{thm:rule-vel-char} recasts the question of whether two rules of the same order have the same trajectories in terms of the quantities $q_{J,\rul}$, so understanding these quantities is key. On the one hand, the presence of indicator functions $\ind{ab \in E(H)}$ in the general expression for $q_{J,\rul}$ given by~\eqref{eq:rule-orbit-coeff} makes it difficult to work with directly. On the other hand, we observe that the symmetries present in symmetric rules allow us to obtain cleaner expressions for the quantities $q_{J,\rul}$; this is the focus of Lemma~\ref{lem:sym-replace-density}. Indeed, we can directly verify the $(\Leftarrow)$ direction of Theorem~\ref{thm:rule-unique-symm-det} by combining Theorem~\ref{thm:rule-vel-char} with the simplified expressions from this key lemma. For the $(\Rightarrow)$ direction, the simplified expressions from Lemma~\ref{lem:sym-replace-density} give us enough leeway and control to construct `perturbed' versions of symmetric non-deterministic rules of order $k\ge3$ which have the same trajectories. To handle non-symmetric rules and complete the proof, we observe and directly verify that every such rule has a natural symmetric version with the same trajectories. This gives our characterization of unique rules.

To prepare for our proof of Theorem~\ref{thm:rule-unique-symm-det}, we introduce some notation, terminology and useful lemmas. The proof of Theorem~\ref{thm:rule-unique-symm-det} is given in Section~\ref{subsec:proof-thm-unique}. We use the language of group actions to describe the symmetries observed. Let a group $\Gamma$ act on a set $X$. For $x\in X$ the \emph{orbit} of $x$ is $\orb_\Gamma(x) = \{\gamma\cdot x : \gamma\in\Gamma\}$ and the \emph{stabilizer} of $x$ is $\stab_\Gamma(x) = \{\gamma\in\Gamma : \gamma\cdot x = x \}$. Set $\ORB_\Gamma(X) = \{\orb_{\Gamma}(x) : x \in X\}$. Let $\Sigma$ be a subgroup of $\Gamma$. A \emph{left coset} of $\Sigma$ in $\Gamma$ is a set $\gamma\Sigma = \{\gamma\sigma : \sigma\in\Sigma\}$ with $\gamma\in\Gamma$. The following technical lemma allows us to recast sums over orbits of group actions as sums over groups.

\begin{lemma} \label{lem:sums-orbit-group}
Given an action of a group $\Gamma$ on a set $X$, a function $f \colon X \to \R$ and an element $x\in X$, we have
\begin{equation*}
\sum_{\gamma\in\Gamma}\frac{f(\gamma\cdot x)}{|\Gamma|} = \sum_{y\in\orb_{\Gamma}(x)}\frac{f(y)}{|\orb_{\Gamma}(x)|}\;.
\end{equation*}
\end{lemma}

\begin{proof}
Observe that the left cosets of $\Sigma := \stab_{\Gamma}(x)$ partition $\Gamma$, the map $\gamma\Sigma \mapsto \gamma\cdot x$ gives a well-defined bijection from the collection $\Gamma/\Sigma$ of left cosets of $\Sigma$ in $\Gamma$ to $\orb_{\Gamma}(x)$ and for each $\gamma\in\Gamma$ the map $\sigma \mapsto \gamma\sigma$ gives a bijection from $\Sigma$ to $\gamma\Sigma$, so by grouping terms according to left cosets of $\Sigma$ we have
\begin{equation*}
\sum_{\gamma\in\Gamma}\frac{f(\gamma\cdot x)}{|\Gamma|} = \sum_{\rho\Sigma\in \Gamma/\Sigma}\sum_{\gamma\in\rho\Sigma}\frac{f(\gamma\cdot x)}{|\Gamma|} = \sum_{y\in\orb_{\Gamma}(x)}\frac{f(y)}{|\orb_{\Gamma}(x)|}
\end{equation*}
as required.
\end{proof}

Now we describe the group actions we work with. In addition to the action of $S_k$ on $\lgr{k}$ described just before Theorem~\ref{thm:rule-unique-symm-det}, we also have the following actions of $S_k$ on $\binom{[k]}{2}$ and $[k]^{(2)}$. For $\sigma \in S_k$ and $ab\in\binom{[k]}{2}$ we set $\sigma \cdot ab = \sigma(a)\sigma(b)$. For $\sigma \in S_k$ and $(a,b)\in[k]^{(2)}$ we set $\sigma \cdot (a,b) = (\sigma(a),\sigma(b))$. For $\ell\in\N$ and $A_1,\dots,A_{\ell}\in\{\lgr{k},\binom{[k]}{2},[k]^{(2)}\}$ we have the following action of $S_k$ on $\A=\prod_{i\in[\ell]}A_i$. For $\sigma \in S_k$ and $(a_i)_{i\in[\ell]}\in\A$ set $\sigma \cdot (a_i)_{i\in[\ell]} = (\sigma \cdot a_i)_{i\in[\ell]}$. Recall that $\cg_k = \lgr{k} \times [k]^{(2)}$. Observe that $\cj_k$ is simply the set $\ORB_{S_k}(\cg_k)$ of orbits for the action of $S_k$ on $\cg_k$. Note that every group action of $S_k$ on a set $X$ naturally induces a group action of each subgroup $\Gamma$ of $S_k$ on the same set $X$.

To prepare for our key lemma, we introduce some notation. Let $F\in\lgr{k}$. Set $S_F = \stab_{S_k}(F)$ and $B_F = \ORB_{S_F}(\binom{[k]}{2})$; the automorphism group of $F$ corresponds to $S_F$. Since we are primarily interested in symmetric rules, it makes sense to think about pairs of vertices in the replacement graph not individually but rather grouped up according to the orbits in $B_F$. In particular, the number of edges on each orbit is a key parameter. Let $\PP_F = \prod_{I\in B_F}[|I|]_0$. Assign to each $H\in\lgr{k}$ a vector $\p^{F,H} = (p^{F,H}_I)_{I\in B_F}\in\PP_F$ given by $p^{F,H}_I = |E(H) \cap I|$ for all $I\in B_F$. For each vector $\p\in\PP_F$ set $\lgr{k}(F,\p) = \{H\in\lgr{k}:\p^{F,H}=\p\}$ and $\rul(F,\p)=\sum_{H\in\lgr{k}(F,\p)} \rul_{F,H}$. For $I\in B_F$ and $0\le\ell\le|I|$ set $\lgr{k}(F,I,\ell) = \{H\in\lgr{k}:\p^{F,H}_I=\ell\}$ and $\rul(F,I,\ell)=\sum_{H\in\lgr{k}(F,I,\ell)}\rul_{F,H}$. For $ab\in\binom{[k]}{2}$ and $0\le\ell\le|\orb_{S_F}(ab)|$ set $\lgr{k}(F,ab,\ell) = \lgr{k}(F,\orb_{S_F}(ab),\ell)$ and $\rul(F,ab,\ell) = \rul(F,\orb_{S_F}(ab),\ell)$. Our key lemma gives useful expressions for $q_{J,\rul}$ when $\rul$ is symmetric.

\begin{lemma} \label{lem:sym-replace-density}
Given a symmetric rule $\rul$ of order $k$, a graph $F\in\lgr{k}$ and $ab\in\binom{[k]}{2}$ we have
\begin{equation*}
\begin{split}
q_{\orb_{S_k}(F^{a,b}),\rul} & = |\orb_{S_k}(F^{a,b})|\left(\sum_{H\in\lgr{k}} \rul_{F,H} \cdot \ind{ab \in E(H)} - \ind{ab \in E(F)} \right) \\
& = |\orb_{S_k}(F^{a,b})|\left(\sum_{\ell=0}^{|\orb_{S_F}(ab)|} \frac{\ell\cdot\rul(F,ab,\ell)}{|\orb_{S_F}(ab)|} - \ind{ab \in E(F)} \right).
\end{split}
\end{equation*}
\end{lemma}

\begin{proof}
For $G^{c,d}\in\orb_{S_k}(F^{a,b})$ pick $\sigma \in S_k$ such that $G=\sigma \cdot F$ and $cd = \sigma \cdot ab$. Clearly we have $\ind{cd \in E(G)} = \ind{ab \in E(F)}$. Since $\rul$ is symmetric and $H\mapsto\sigma \cdot H$ defines a bijection on $\lgr{k}$, we may replace sums over $H\in\lgr{k}$ with sums over $\sigma \cdot H\in\lgr{k}$ and obtain
\[\sum_{H\in\lgr{k}} \rul_{G,H} \cdot \ind{cd \in E(H)} = \sum_{H\in\lgr{k}} \rul_{G,\sigma \cdot H} \cdot \ind{\sigma \cdot ab \in E(\sigma \cdot H)} = \sum_{H\in\lgr{k}} \rul_{F,H} \cdot \ind{ab \in E(H)}.\]
Hence, by~\eqref{eq:rule-orbit-coeff} we have
\begin{equation} \label{eq:rule-coeff-sym}
q_{\orb_{S_k}(F^{a,b}),\rul} = |\orb_{S_k}(F^{a,b})|\left(\sum_{H\in\lgr{k}} \rul_{F,H} \cdot \ind{ab \in E(H)} - \ind{ab \in E(F)} \right)\;.
\end{equation}

Take $G \in J \in \ORB_{S_F}(\lgr{k})$. By replacing sums over $\sigma \in S_F$ with sums over $\sigma^{-1} \in S_F$ and applying Lemma~\ref{lem:sums-orbit-group} for $G\in\lgr{k}$ and for $ab\in\binom{[k]}{2}$, we obtain
\begin{equation*}
\sum_{H\in\orb_{S_F}(G)}\frac{\ind{ab \in E(H)}}{|\orb_{S_F}(G)|} = \sum_{\sigma\in S_F} \frac{\ind{ab \in E(\sigma \cdot G)}}{|S_F|} = \sum_{\sigma\in S_F} \frac{\ind{\sigma \cdot ab \in E(G)}}{|S_F|} = \sum_{cd\in\orb_{S_F}(ab)}\frac{\ind{cd \in E(G)}}{|\orb_{S_F}(ab)|}\;.
\end{equation*}
Rearranging yields
\begin{equation} \label{eq:perm-orbit-edges}
\sum_{H\in\orb_{S_F}(G)}\ind{ab \in E(H)} = \frac{|\orb_{S_F}(G)|}{|\orb_{S_F}(ab)|}\cdot\sum_{cd\in\orb_{S_F}(ab)}\ind{cd \in E(G)}\;.
\end{equation}
Now $\rul$ is symmetric, so for all $G\in\lgr{k}$ and $H\in\orb_{S_F}(G)$ we have $\rul_{F,G} = \rul_{F,H}$. Hence, in combination with~\eqref{eq:rule-coeff-sym} and~\eqref{eq:perm-orbit-edges} we obtain
\[ q_{\orb_{S_k}(F^{a,b}),\rul} = |\orb_{S_k}(F^{a,b})|\left(\sum_{\ell=0}^{|\orb_{S_F}(ab)|} \frac{\ell\cdot\rul(F,ab,\ell)}{|\orb_{S_F}(ab)|} - \ind{ab \in E(F)} \right) \]
as desired.
\end{proof}

\subsection{Proof of Theorem~\ref{thm:rule-unique-symm-det}} \label{subsec:proof-thm-unique}

We present each direction of the proof of Theorem~\ref{thm:rule-unique-symm-det} separately.

\begin{proof}[Proof of Theorem~\ref{thm:rule-unique-symm-det} $(\Leftarrow)$]
We first consider when $\rul$ is a rule of order $2$. Let $\rul'$ be a rule of the same order with the same trajectories as $\rul$. Take $F\in\lgr{2}$ and set $J = \orb_{S_2}(F^{1,2}) = \{F^{1,2},F^{2,1}\}$. By~\ref{item:same-coeff} of Theorem~\ref{thm:rule-vel-char} we have 
\[2\left(\rul'_{F,K_2} - \ind{12 \in E(F)}\right) = q_{J,\rul'} = q_{J,\rul} = 2\left(\rul_{F,K_2} - \ind{12 \in E(F)}\right)\;,\]
so we obtain $\rul'_{F,K_2} = \rul_{F,K_2}$. But $\lgr{2} = \{K_2,\overline{K_2}\}$, so this gives $\rul = \rul'$. Hence, $\rul$ is unique.

Next, we consider when $\rul$ is a symmetric deterministic rule of order $k$. Let $\rul'$ be a rule of the same order with the same trajectories as $\rul$. Let $F\in\lgr{k}$ and fix the unique $H_0\in\lgr{k}$ such that $\rul_{F,H_0}=1$. Take $ab\in\binom{[k]}{2}$ and set $J=\orb_{S_k}(F^{a,b})$. By~\ref{item:same-coeff} of Theorem~\ref{thm:rule-vel-char} we have $q_{J,\rul'} = q_{J,\rul}$, so by Lemma~\ref{lem:sym-replace-density} we have $\sum_{H\in\lgr{k}} \rul'_{F,H} \cdot (\ind{ab \in E(H)} - \ind{ab \in E(H_0)}) = 0$. It follows that $\rul'_{F,H} = 0$ for all $H\neq H_0$. This gives $\rul = \rul'$, so $\rul$ is unique.
\end{proof}

\begin{proof}[Proof of Theorem~\ref{thm:rule-unique-symm-det} $(\Rightarrow)$]
We prove the contrapositive statement. We first consider when $\rul$ is not symmetric. Set $\ca = \lgr{k} \times \lgr{k}$ and $\cb = \ORB_{S_k}(\ca)$. The \emph{symmetric version} of $\rul$ is the matrix $\rul^{\times} = (\rul^{\times}_{F,H})_{F,H\in\lgr{k}}$ where for all $(F,H) \in B \in \cb$ we have $\rul^{\times}_{F,H} = \frac{1}{|B|}\sum_{(F',H') \in B}\rul_{F',H'}$. The following claim shows that $\rul^{\times}$ is a symmetric rule of the same order with the same trajectories as $\rul$; since $\rul$ is not symmetric, this shows that $\rul$ is not unique.

\begin{claim} \label{claim:sym-ver-rule-traj}
$\rul^{\times}$ is a symmetric rule with the same trajectories as $\rul$.
\end{claim}

\begin{claimproof}
We first verify that $\rul^{\times}$ is a symmetric rule. By the definition of $\rul^{\times}$ it is symmetric and we have $\rul^{\times}_{F,H}\in[0,1]$ for all $(F,H)\in\ca$, so it remains to show that for all $F\in\lgr{k}$ we have $\sum_{H\in\lgr{k}}\rul^{\times}_{F,H} = 1$.

Fix $P\in\ORB_{S_k}(\lgr{k})$. Note that $\cc = \{A\in\ORB_{S_k}(\ca):A\cap(P\times\lgr{k})\ne\nth\}$ partitions $P\times\lgr{k}$. Indeed, this follows from the fact that the collection $\ORB_{S_k}(\ca)$ of orbits partitions $\ca$ and for all $(G,H) \in P\times\lgr{k}$ we have $\orb_{S_k}(G,H)\subseteq P\times\lgr{k}$. Hence, by the definition of $\rul^{\times}$ we have
\begin{equation} \label{eq:rule-orbit-row-sum}
\sum_{(G,H)\in P\times\lgr{k}}\rul^{\times}_{G,H} = \sum_{A\in\cc} \sum_{(G,H)\in A}\rul^{\times}_{G,H} = \sum_{A\in\cc} \sum_{(G,H)\in A}\rul_{G,H} = \sum_{(G,H)\in P\times\lgr{k}}\rul_{G,H} = |P|\;.
\end{equation}
Now take $F\in P$. By the definition of $P$, for each $F'\in P$ we may pick $\sigma \in S_k$ such that $F' = \sigma \cdot F$, so we have $\sum_{H\in\lgr{k}}\rul^{\times}_{F',H} = \sum_{H\in\lgr{k}}\rul^{\times}_{\sigma \cdot F,\sigma \cdot H} = \sum_{H\in\lgr{k}}\rul^{\times}_{F,H}$. In combination with~\eqref{eq:rule-orbit-row-sum}, we obtain $|P| = \sum_{(G,H)\in P\times\lgr{k}}\rul^{\times}_{G,H} = |P| \cdot \sum_{H\in\lgr{k}}\rul^{\times}_{F,H}$, which implies $\sum_{H\in\lgr{k}}\rul^{\times}_{F,H} = 1$. Hence, $\rul^{\times}$ is a rule as required.

Now we show that $\rul^{\times}$ and $\rul$ have the same trajectories by verifying the conditions of~\ref{item:same-coeff}. Let $J\in\cj_k$ and $(F^{a,b},H)\in J\times\lgr{k}$. Set $A := \orb_{S_k}(F^{a,b},H)$ and $B := \orb_{S_k}(F,H)$. Noting that $\ind{\sigma \cdot ab \in E(\sigma \cdot H)} = \ind{ab \in E(H)}$ for all $\sigma \in S_k$ and applying Lemma~\ref{lem:sums-orbit-group} for $(F^{a,b},H) \in \cg_k\times\lgr{k}$ and for $(F,H) \in \lgr{k}\times\lgr{k}$, we obtain
\begin{equation*}
\begin{split}
\sum_{(G^{c,d},H')\in A} \frac{\rul_{G,H'}\ind{cd \in E(H')}}{|A|} & = \sum_{\sigma\in S_k} \frac{\rul_{\sigma \cdot F,\sigma \cdot H}\ind{\sigma \cdot ab \in E(\sigma \cdot H)}}{|S_k|} = \sum_{\sigma\in S_k} \frac{\rul_{\sigma \cdot F,\sigma \cdot H}\ind{ab \in E(H)}}{|S_k|} \\
& = \sum_{(G,H')\in B} \frac{\ind{ab \in E(H)}\rul_{G,H'}}{|B|}\;.
\end{split}
\end{equation*}
Since the analogous equation for $\rul^{\times}$ also holds and by the definition of $\rul^{\times}$ we have
\begin{equation*}
\sum_{(F,H)\in B}\rul_{F,H} = \sum_{(F,H)\in B}\rul^{\times}_{F,H}\;,
\end{equation*}
we obtain
\begin{equation*}
\sum_{(F^{a,b},H)\in A}\rul_{F,H}\ind{ab \in E(H)} = \sum_{(F^{a,b},H)\in A}\rul^{\times}_{F,H}\ind{ab \in E(H)}\;.
\end{equation*}
Hence, by~\eqref{eq:rule-orbit-coeff} we have
\begin{equation*}
\begin{split}
q_{J,\rul} & = \sum_{A\in\ORB_{S_k}(J\times\lgr{k})} \sum_{(F^{a,b},H) \in A} \rul_{F,H} \cdot \ind{ab \in E(H)} - \sum_{F^{a,b}\in J}\ind{ab \in E(F)} \\
& = \sum_{A\in\ORB_{S_k}(J\times\lgr{k})} \sum_{(F^{a,b},H) \in A} \rul^{\times}_{F,H} \cdot \ind{ab \in E(H)} - \sum_{F^{a,b}\in J}\ind{ab \in E(F)} = q_{J,\rul^{\times}}.
\end{split}
\end{equation*}
This establishes that the conditions of~\ref{item:same-coeff} hold, so by Theorem~\ref{thm:rule-vel-char} the rules $\rul$ and $\rul^{\times}$ have the same trajectories as required.
\end{claimproof}

The only rule of order $1$ is symmetric and deterministic, so it remains to consider when $\rul$ is a non-deterministic symmetric rule of order $k\ge3$. First we prove the following claim.

\begin{claim} \label{claim:detsymrule-orbit-fullempty}
A symmetric rule $\rul$ is deterministic if and only if for all $F\in\lgr{k}$ and $I\in B_F$ there exists $q_{F,I}\in\{0,|I|\}$ such that for all $H\in\lgr{k}$ such that $p^{F,H}_I\ne q_{F,I}$ we have $\rul_{F,H}=0$.
\end{claim}

\begin{claimproof}
We start with the $(\Rightarrow)$ direction. Since $\rul$ is deterministic, we have a unique $G_F\in\lgr{k}$ such that $\rul_{F,H} = \ind{H=G_F}$. Take $F\in\lgr{k}$ and $I\in B_F$. We set $q_{F,I}:=|E(G_F)\cap I|$ and claim that $q_{F,I}\in\{0,|I|\}$. Indeed, suppose that $q_{F,I}>0$. Then we may pick $e\in E(G_F)\cap I$. Now take $\sigma\in S_F$. Since $\rul$ is symmetric, we have $1 = \rul_{F,G_F} = \rul_{\sigma\cdot F,\sigma\cdot G_F} = \rul_{F,\sigma\cdot G_F}$, which implies $G_F = \sigma\cdot G_F$. In particular, we have $\sigma\cdot e\in E(G_F)$. Since $I = \orb_{S_F}(e)$, we have $E(G_F)\cap I = I$ as required. Now take $H\in\lgr{k}$ such that $p^{F,H}_I=|E(H)\cap I|\ne |E(G_F)\cap I| = q_{F,I}$. In particular, we have $H\ne G_F$. Since $\rul$ is deterministic, we have $\rul_{F,H}=0$.

Now we consider the $(\Leftarrow)$ direction. Fix $F\in\lgr{k}$ and for all $I\in B_F$ take $q_{F,I}\in\{0,|I|\}$ whose existence is guaranteed by assumption. Let $G_F\in\lgr{k}$ be the labelled graph on $[k]$ such that $E(G_F) = \bigcup_{I:q_{F,I}=|I|}I$. Now $G_F$ is the unique labelled graph $H\in\lgr{k}$ such that $p^{F,H}_I = |E(H)\cap I| = q_{F,I}$ for all $I\in B_F$, so for all $H\in\lgr{k}$ we have $\rul_{F,H} = \ind{H=G_F}$. It follows that $\rul$ is deterministic.
\end{claimproof}

It follows from Claim~\ref{claim:detsymrule-orbit-fullempty} that the following (not necessarily mutually exclusive) cases cover all non-deterministic symmetric rules of order $k\ge3$.
\begin{enumerate}[label=(D\arabic{*})]
\item \label{item:rule-intr-supp} There are $F\in\lgr{k}$, $I\in B_F$, $\p\in\PP_F$ and $H\in\lgr{k}(F,\p)$ such that $|I|\ge2$, $p_I\in[|I|-1]$ and $\rul_{F,H}>0$.
\item \label{item:rule-endpts-supp} There are $F\in\lgr{k}$, $I\in B_F$, $\p^{(i)}\in\PP_F$ and $H^{(i)}\in\lgr{k}(F,\p^{(i)})$ for all $i\in\{0,1\}$ such that $|I|\ge2$ and the following hold for all $i\in\{0,1\}$. We have $p^{(i)}_I = i|I|$, $p^{(i)}_A=p^{(0)}_A$ for all $A\ne I$ and $\rul_{F,H^{(i)}}>0$.
\item \label{item:rule-bin-supp-symm} There are $F\in\lgr{k}$, $\{ab\}\in B_F$, $\p^{-},\p^{+}\in\PP_F$, $H^{-}\in\lgr{k}(F,\p^{-})$ and $H^{+}\in\lgr{k}(F,\p^{+})$ such that $p^{-}_{\{ab\}}=0$, $p^{+}_{\{ab\}}=1$, $p^{-}_A=p^{+}_A$ for $A\ne\{ab\}$ and $\rul_{F,H^{-}},\rul_{F,H^{+}}>0$.
\end{enumerate}
Since $\lgr{k}(F,\p)$ is the set of $H\in\lgr{k}$ satisfying $|E(H)\cap I|=p_I$ for all $I\in B_F$, it is nonempty for any $\p$ satisfying $0\le p_I\le|I|$ for all $I\in B_F$.

For~\ref{item:rule-intr-supp} define $\p^{-},\p^{+}\in\PP_F$ by setting $p^{-}_I=p_I-1$, $p^{+}_I=p_I+1$ and $p^{-}_A,p^{+}_A=p_A$ for all $A\ne I$. Pick $H^{-}\in\lgr{k}(F,\p^{-})$ and $H^{+}\in\lgr{k}(F,\p^{+})$. Take $0<\eps\le|\orb_{S_k}(F,H)|\rul_{F,H}$ and define a modification $\rul^{\eps}$ of $\rul$ as follows. For $(F',H')\in\lgr{k}\times\lgr{k}$ we set
\begin{equation*}
\rul^{\eps}_{F',H'}=
\begin{cases}
\rul_{F',H'} - \frac{\eps}{|\orb_{S_k}(F,H)|} & \textrm{if }(F',H')\in\orb_{S_k}(F,H) \\
\rul_{F',H'} + \frac{\eps}{2|\orb_{S_k}(F,H^{-})|} & \textrm{if }(F',H')\in\orb_{S_k}(F,H^{-}) \\
\rul_{F',H'} + \frac{\eps}{2|\orb_{S_k}(F,H^{+})|} & \textrm{if }(F',H')\in\orb_{S_k}(F,H^{+}) \\
\rul_{F',H'} & \textrm{otherwise.}
\end{cases}
\end{equation*}
It is straightforward to check that $\rul^{\eps}$ is a symmetric rule distinct from $\rul$. To show that $\rul$ is not unique, we shall show that $\rul$ and $\rul^{\eps}$ have the same trajectories by verifying the conditions of~\ref{item:same-coeff} from Theorem~\ref{thm:rule-vel-char}.

We consider three categories of $J\in\cj_k$. First, we have $J\ne\orb_{S_k}(F^{a,b})$ for all $ab \in \binom{[k]}{2}$. In this case, by~\eqref{eq:rule-orbit-coeff} and the definition of $\rul^{\eps}$ we have $q_{J,\rul} = q_{J,\rul^{\eps}}$. Next, we have $J=\orb_{S_k}(F^{a,b})$ for some $ab \notin I$. Here, by the definition of $\rul^{\eps}$ we have $\rul^{\eps}(F,ab,\ell) = \rul(F,ab,\ell)$ for all $\ell$, so by Lemma~\ref{lem:sym-replace-density} we have $q_{J,\rul} = q_{J,\rul^{\eps}}$. Finally, we have $J=\orb_{S_k}(F^{a,b})$ for some $ab \in I$. In this case, we have $\rul^{\eps}(F,ab,p_I) = \rul(F,ab,p_I) - \eps$, $\rul^{\eps}(F,ab,p_I-1) = \rul(F,ab,p_I-1) + \frac{\eps}{2}$, $\rul^{\eps}(F,ab,p_I+1) = \rul(F,ab,p_I+1) + \frac{\eps}{2}$ and $\rul^{\eps}(F,ab,\ell) = \rul(F,ab,\ell)$ otherwise. Hence, by Lemma~\ref{lem:sym-replace-density} we have $q_{J,\rul} = q_{J,\rul^{\eps}}$. Since the conditions of~\ref{item:same-coeff} from Theorem~\ref{thm:rule-vel-char} hold, $\rul$ and $\rul^{\eps}$ have the same trajectories and $\rul$ is not unique.

For~\ref{item:rule-endpts-supp} take $0<\eps\le\min\{|\orb_{S_k}(F,H^{(0)})|\rul_{F,H^{(0)}},|\orb_{S_k}(F,H^{(1)})|\rul_{F,H^{(1)}}\}$ and distinguish two subcases: $|I|\ge3$ and $|I|=2$. We begin with the former. Define $\q^{-},\q^{+}\in\PP_F$ by setting $q^{-}_I=1$, $q^{+}_I=|I|-1$ and $q^{-}_A=q^{+}_A=p^{(0)}_A$ for all $A\ne I$. Pick $G^{-}\in\lgr{k}(F,\q^{-})$ and $G^{+}\in\lgr{k}(F,\q^{+})$. Define a modification $\rul^{\eps}$ of $\rul$ as follows. For $(F',H')\in\lgr{k}\times\lgr{k}$ we set
\begin{equation*}
\rul^{\eps}_{F',H'}=
\begin{cases}
\rul_{F',H'} -\frac{\eps}{|\orb_{S_k}(F,H^{(i)})|} & \textrm{if }(F',H')\in\orb_{S_k}(F,H^{(i)}) \\
\rul_{F',H'} + \frac{\eps}{|\orb_{S_k}(F,G^{-})|} & \textrm{if }(F',H')\in\orb_{S_k}(F,G^{-}) \\
\rul_{F',H'} + \frac{\eps}{|\orb_{S_k}(F,G^{+})|} & \textrm{if }(F',H')\in\orb_{S_k}(F,G^{+}) \\
\rul_{F',H'} & \textrm{otherwise.}
\end{cases}
\end{equation*}
In a manner similar to that for~\ref{item:rule-intr-supp}, we may verify that $\rul^{\eps}$ is a symmetric rule distinct from $\rul$ which has the same trajectories, thereby establishing that $\rul$ is not unique.

For the latter subcase define $\q\in\PP_F$ by setting $q_I=1$ and $q_A=p^{(0)}_A$ for all $A\ne I$. Pick $G\in\lgr{k}(F,\q)$ and define a modification $\rul^{\eps}$ of $\rul$ as follows. For $(F',H')\in\lgr{k}\times\lgr{k}$ we set
\begin{equation*}
\rul^{\eps}_{F',H'}=
\begin{cases}
\rul_{F',H'} -\frac{\eps}{|\orb_{S_k}(F,H^{(i)})|} & \textrm{if }(F',H')\in\orb_{S_k}(F,H^{(i)}) \\
\rul_{F',H'} + \frac{2\eps}{|\orb_{S_k}(F,G)|} & \textrm{if }(F',H')\in\orb_{S_k}(F,G) \\
\rul_{F',H'} & \textrm{otherwise.}
\end{cases}
\end{equation*}
In a manner similar to that for~\ref{item:rule-intr-supp}, we may verify that $\rul^{\eps}$ is a symmetric rule distinct from $\rul$ which has the same trajectories, thereby establishing that $\rul$ is not unique.

For~\ref{item:rule-bin-supp-symm} pick $G\in\lgr{k}$ distinct from $F$ so that for some $cd\in\binom{[k]}{2}$ we have $G^{c,d}\in\orb_{S_k}(F^{a,b})$; we may always do so because $k\ge3$ and $\orb_{S_F}(ab) = \{ab\}$. Fix $\sigma\in S_k$ such that $G = \sigma\cdot F$, take $0<\eps\le\min\{\rul_{F,H^{-}},\rul_{G,H^{+}},1-\rul_{F,H^{+}},1-\rul_{G,H^{-}}\}$. Define a modification $\rul^{\eps}$ of $\rul$ as follows. For $(F',H')\in\lgr{k}\times\lgr{k}$ we set
\begin{equation*}
\rul^{\eps}_{F',H'}=
\begin{cases}
\rul_{F',H'} - \eps & \textrm{if }(F',H')=(F,H^{-}),(G,\sigma\cdot H^{+}) \\
\rul_{F',H'} + \eps & \textrm{if }(F',H')=(F,H^{+}),(G,\sigma\cdot H^{-}) \\
\rul_{F',H'} & \textrm{otherwise.}
\end{cases}
\end{equation*}
In a manner similar to that for~\ref{item:rule-intr-supp}, we may verify that $\rul^{\eps}$ is a rule distinct from $\rul$ which has the same trajectories, thereby establishing that $\rul$ is not unique.
\end{proof}

\section{Concluding remarks}

In this paper, we obtain a complete characterization of the equivalence classes of flip process rules with the same graphon trajectories and characterize the flip process rules which are unique in their equivalence classes. In this section, we shall discuss a number of possible directions for further study.

\subsection{Characterization by flip process distributions} \label{subsec:flip-char}

It is natural to return to the flip process setting and consider the flip process distribution analogue of Question~\ref{qn:rule-traj-equal}.

\begin{question} \label{qn:rule-flip-equal}
When do two rules $\rul_1$ and $\rul_2$ have the same flip process distributions?
\end{question}

Since flip processes are time-homogeneous Markov processes, it suffices to consider their one-step evolutions, that is, the distribution of $G_1$ given an initial graph $G_0 = G$. This is analogous to studying the velocity operator for graphon trajectories and leads us to the following conjecture. Together with the `lifting' procedure relating rules of different orders described just below Theorem~\ref{thm:rule-vel-char-intro}, a positive answer would fully resolve Question~\ref{qn:rule-flip-equal}.

\begin{conj} \label{conj:flip-char}
The following are equivalent for rules $\rul$ and $\rul'$ of the same order $k\in\N$.
\begin{enumerate}[label=(K\arabic{*})]
\item \label{item:same-symm} For all $B\in\ORB_{S_k}(\lgr{k} \times \lgr{k})$ we have $\sum_{(F',H') \in B}\rul_{F',H'} = \sum_{(F',H') \in B}\rul'_{F',H'}$\;.
\item \label{item:same-flip-dist} $\rul$ and $\rul'$ have the same flip process distributions.
\end{enumerate}
\end{conj}

One potential approach to resolving Conjecture~\ref{conj:flip-char} is to replicate the proof framework of Theorem~\ref{thm:rule-vel-char}. Indeed, much of the proof has a graph theoretic flavour and would carry over by analogy. In particular, one could prove an analogue of Lemma~\ref{lem:rooted-density-quantify} for a suitable analogue of the operator $\Troot{F^{a,b}}$. However, the flip process sampling method means that we would need to count \emph{injective} copies of graphs and so we would get terms $\prod_{i\in[m]}(z-i+1)$ resembling falling factorials instead of simple monomials $z^m$. As a consequence, showing that we have a nonzero polynomial function would require more work.

\subsection{Trajectories with varying speeds}

Our definition of rules $\rul_1$ and $\rul_2$ with the same trajectories requires $\trajrul{\rul_1}{t}{W} = \trajrul{\rul_2}{t}{W}$ for all $t$. That is, the trajectories starting at $W$ have to be the same in `real-time' under both $\rul_1$ and $\rul_2$. One might ask what happens if we were to relax this time-restriction and allow speed variation. Here we state an example question. 

\begin{question} \label{qn:rule-traj-equal-time-fn}
When do two rules $\rul_1$ and $\rul_2$ satisfy $\trajrul{\rul_1}{t}{W} = \trajrul{\rul_2}{f(W,t)}{W}$ for all $W\in\Gra$ and $t\in[0,\infty)$ for a suitable `time function' $f\colon\Gra\times[0,\infty)\to[0,\infty)$?
\end{question}

One notable case is that of uniform time dilation, that is, when $f(W,t) = Ct$ for some constant $C>0$. Here the proof framework of Theorem~\ref{thm:rule-vel-char} can be adapted to give the characterization $q_{J,\rul_1} = C \cdot q_{J,\rul_2}$ (analogous to~\ref{item:same-coeff}). There are concrete examples of this form of speed variation. Indeed, let $\rul_1$ be the triangle removal rule and $\rul_2$ be the rule which deletes a random edge from any sampled triangle (and does nothing otherwise). It is straightforward to check that $\trajrul{\rul_1}{t}{W} = \trajrul{\rul_2}{3t}{W}$.

We remark that there are also concrete examples of non-uniform time dilation. Let $\rul_1$ and $\rul_2$ be rules of order $6$ which act as follows. If both sets $\{v_1,v_2,v_3\}$ and $\{v_4,v_5,v_6\}$ of drawn vertices induce a triangle, $\rul_1$ deletes the edges of the triangle on $\{v_1,v_2,v_3\}$. If the drawn vertices $\{v_1,v_2,v_3\}$ induce a triangle, $\rul_2$ deletes the edges of that triangle. Otherwise, they do nothing. It is straightforward to check that $\vel[\rul_1]W = f(W) \cdot \vel[\rul_2]W$, where $f(W)$ represents the triangle density of $W$. By Definition~\ref{def:traj} and Theorem~\ref{thm:traj}, the trajectories $\trajrul{\rul_1}{t}{W}$ and $\trajrul{\rul_2}{t}{W}$ are unique and satisfy $\trajrul{\rul_1}{t}{W} = \trajrul{\rul_2}{g(t,W)}{W}$ with $g(t,W) = \int_0^tf(\trajrul{\rul_1}{s}{W}) \D s$.

\section*{Acknowledgements}

The author would like to thank Jan Hladký for helpful discussions and insightful comments.

\bibliographystyle{siam}
\bibliography{flip-rule-bib}

\end{document}